\documentclass[review,12pt]{elsarticle}

\usepackage{amssymb}

\usepackage{xspace}
\usepackage{algorithm}
\usepackage{algorithmic}
\usepackage{subfigure}
\usepackage{amsfonts}
\usepackage{amsmath}
\usepackage{enumerate}
\usepackage{framed}
\usepackage{graphicx}
\usepackage{lscape}
\usepackage{bm}
\usepackage{color}
\usepackage{multicol}
\usepackage{amsthm}

\newcommand{\Leja}{{L\'{e}ja}\xspace}

\newcommand{\algref}[1]{{Algorithm~\ref{#1}}}
\newcommand{\figref}[1]{{Figure~\ref{#1}}}
\newcommand{\tabref}[1]{{Table~\ref{#1}}}

\renewcommand{\eqref}[1]{{(\ref{#1})}}

\begin{document}

\begin{frontmatter}



\title{Efficient Simulation of Geothermal Processes in Heterogeneous Porous Media based on the Exponential Rosenbrock-Euler  and Rosenbrock-type Methods}
 


\author[at]{A.~Tambue\corref{cor1}}
\ead{Antoine.Tambue@math.uib.no}
\cortext[cor1]{Corresponding author}
\address[at]{ Department of Mathematics, University of Bergen, P.O. Box 7800, N-5020 Bergen, Norway}
\author[at,ib]{I.~Berre}
\ead{Inga.Berre@math.uib.no}
\address[ib]{Christian Michelsen Research AS, Norway}
\author[at,jmn]{J.M.Nordbotten}
\ead{Jan.Nordbotten@math.uib.no}
\address[jmn]{Department of Civil and Environmental Engineering, Princeton University, USA}

\address{}

\begin{abstract}
Simulation of geothermal systems is challenging due to coupled physical processes in highly
heterogeneous media. Combining the exponential  Rosenbrock--Euler
and Rosenbrock-type methods with control-volume (two-point flux approximation) space
discretizations leads to efficient numerical techniques for simulating geothermal systems. In terms of efficiency
and accuracy, the exponential Rosenbrock--Euler  time integrator has advantages over standard time-dicretization
schemes, which suffer from time-step restrictions or excessive numerical diffusion when advection
processes are dominating. Based on linearization of the equation at each time step, we make use of 
matrix exponentials of the Jacobian from the spatial discretization, which %
provide the exact solution in time for the
linearized equations. This is at the expense of
computing the  matrix exponentials of the stiff
Jacobian matrix,
together with propagating a linearized system.
However, using a Krylov subspace or Leja points techniques make
these computations efficient.

The Rosenbrock-type methods  use the appropriate rational functions of the Jacobian 
from the spatial discretization. The parameters in these schemes
are found in consistency with the  required order of convergence in time.
As a result, these schemes are A-stable 
and only a few linear systems are solved at each time step.
 The efficiency of the methods compared
to standard time-discretization techniques are
demonstrated in numerical examples.
\end{abstract}



\end{frontmatter}


\section{Introduction}
\label{intro}

In the subsurface, producing geothermal systems are characterized by coupled hydraulic,
thermal, chemical and mechanical processes. To determine the potential of a geothermal site, and
decide optimal production strategies, it is important to understand and quantify these processes. Rigorous
mathematical modeling and accurate numerical simulations are essential, but multiple
interacting processes acting on different scales leads to challenges in solving the coupled system of
equations.

Standard discretization methods include finite element, control-volume and finite difference methods
for the space discretization 
(see \cite{ MPFA, BundschuhArriaga,ChenHuanMa, FV} and references therein), while standard
implicit, explicit or implicit-explicit methods have until recently mostly been used 
for the discretization in time 
\cite{PodgorneyEtAl2010,Pruess1991}. Challenges with the
discretization are, amongst others, related to severe time-step restrictions associated with explicit
methods and excessive numerical diffusion for implicit methods.
Furthermore, implicit methods require at each time step the solution of 
large systems of nonlinear equations, which may lead to bottlenecks in practical computations.


In this paper, we consider a different approach for the
temporal discretization based on the Exponential Rosenbrock--Euler  Method (EREM) and Rosenbrock-type Methods (ROSM). 
Exponential integrators have recently been suggested as efficient and robust alternatives for the
temporal discretization for several applications
(see \cite{carr,TLGspe1,ATthesis,TLGspe}. Rosenbrock-type methods have been intensively developed in the literature
and used in a variety of applications
(see \cite{Rosen,hairer,langb,lang} and reference therein). 
However, neither of these approaches have yet found wide-spread use in porous media applications.

The mathematical model discussed, consists of a system of partial differential equations that express
conservation of mass and energy. In addition, the model entails phenomenological
laws describing processes active in the reservoir, such as Darcy's law for fluid flow with
variable density and viscosity, Fourier's law of heat conduction, and those describing the relation between
fluid properties (nonlinear fluid expansivity and compressibility) and porosity subject to pressure and
temperature variations. The resulting system of equations is  nonlinear and coupled and
requires sophisticated numerical techniques.

Our solution technique is based on a sequential approach, which decouples the mathematical model. 
An advantage of this approach is that it allows for specialized solvers for unknowns
with different characteristics. As the linearized fully coupled matrices are
often very poorly conditioned such that small time-steps are required, a carefully chosen sequential approach leads to
higher efficiency and accuracy than a simultaneous solution approach 
if the couplings are not too strong \cite{ChenHuanMa}.
 A finite volume method is applied for the space discretization while
 the Exponential Rosenbrock--Euler Method and Rosenbrock-type methods are applied 
to integrate the systems in time based on successive linearizations.  
 
The  Exponential  Rosenbrock--Euler method is based on the linearization of  
 the ODEs resulting from the space discretization at each time step. The linear part is solved exactly in 
time up to a given tolerance in the computation of a matrix exponential function 
of the Jacobian. The nonlinear part is approximated using low-order Taylor expansions.
As in all exponential integrator schemes, the expense is the computation  
of the matrix exponentials of the stiff Jacobian matrix resulting from the 
spatial discretization. Computing matrix exponentials of stiff matrices
 is a notorious problem in numerical analysis
 \cite{CMCVL}, but new developments for both L\'{e}ja points and
Krylov subspace techniques \cite{LE,kry,SID,ATthesis,Antoine} have led
to efficient numerical approaches; see e.g. 
  \cite{calarioster,carr,TLGspe1,ATthesis} and references
therein. Besides, the method is L-stable, and performs well for super-stiff ODEs.

The Rosenbrock-type methods use the appropriate rational functions of the Jacobian 
from the spatial discretization. The parameters in these schemes
are found in consistency with the  required order of convergence in time.
As a result these schemes are A-stable (as will be discussed below)
and only few linear systems are solved at each time step.

The paper is organized as follows. We present the model equations in
Section 2, and the finite volume method for spatial
discretization in Section 3, along with  temporal discretization schemes.
 The implementation of the Exponential Rosenbrock-Euler method is discussed in Section 4. In Section 5 we present
some numerical examples, which also include simulations for a fractured reservoir, and show comparisons to standard approaches, before we draw 
conclusions in Section 6.
\section{Model equations}
\label{sec2}
We assume single-phase flow of water, which allows the energy equation 
to be written in terms of temperature. The model equations are given by 
\begin{eqnarray}
\label{heat}
\left\lbrace \begin{array}{l}
 (1-\phi) \rho_{s} c_{ps}\dfrac{\partial T_{s}}{\partial t}=(1-\phi)\nabla \cdot(\mathbf{k}_{s} \nabla 
T_{s})+(1-\phi)q_{s} +he(T_{f}-T_{s)}\\
\phi \rho_{f} c_{pf}\dfrac{\partial T_{f}}{\partial t}=\phi\nabla \cdot(\mathbf{k}_{f} \nabla 
T_{f}) -\nabla \cdot(\rho_{f} c_{pf} \mathbf{v}T_{f}) +\phi q_{f} +he(T_{s}-T_{f)}\\ 
\end{array}\right.
\end{eqnarray}
(see \cite{coumoud,nielb}). The model equations are given in the bounded spatial domain
 $ \Omega \subset \mathbb{R}^{d},\;d={2,3}$, with boundary $ \partial \Omega$, and in the time interval is $\left[ 0,\tau\right]$.
Here $ \phi $ is the porosity; $he$ is the heat transfer coefficient; $q$ is the heat production;
$\rho$ is the density; $c_{p}$ stand for  heat capacity; $T$ is the temperature;
$\mathbf{k}$  is the thermal conductivity tensor, with the 
subscripts $f$ and $s$ referring to fluid and rock; and $\mathbf{v}$ is the Darcy velocity given by
\begin{eqnarray}
\label{darcy}
 \mathbf{v} =-\dfrac{\mathbf{K}}{\mu} \left(\nabla p-\rho_{f}\mathbf{g}\right),
\end{eqnarray}
where $\mathbf{K}$ is the permeability tensor, $\mu$ is the viscosity,  $\mathbf{g}$ 
is the gravitational acceleration and $p$ the pressure. The mass balance equation for a single-phase fluid is given by
\begin{eqnarray}
\label{masscon}
 \dfrac{\partial \phi \rho_{f}}{\partial t}=-\nabla \cdot( \mathbf{v}\rho_{f}) +Q_{f},
\end{eqnarray}
where $Q [\text{kg/s}]$ is contribution from a source or sink per time unit.
 Assuming in equation \eqref{masscon} that the rock is slightly compressible, the porosity is a function of pressure
and can be expressed as a linear function, yielding
\begin{eqnarray}
\label{masscon2}
\rho_{f}\dfrac{\partial \phi}{\partial t} +\phi \dfrac{\partial \rho_{f}}{\partial t}&=&\nabla \cdot \left( \dfrac{\rho_{f}\mathbf{K}}{\mu} \left(\nabla p-\rho_{f}\mathbf{g}\right)\right)+Q_{f},
\end{eqnarray}
 with 
\begin{eqnarray}
\label{rockcomp}
 \phi=\phi_{0} \left(1+ \alpha_{b} (p-p_{0}) \right),
\end{eqnarray}
where $\phi_{0}$ is the porosity at the initial pressure, $p_{0}$ the initial 
pressure and $\alpha_{b}$  the bulk vertical compressibility of the porous medium.

Notice that
\begin{eqnarray}
\label{masscon3}
 \dfrac{\partial \rho_{f}}{\partial t}&=&\dfrac{\partial \rho_{f}}{\partial T_{f}}\dfrac{\partial T_{f}}{\partial t}+\dfrac{\partial \rho_{f}}{\partial p}\dfrac{\partial p}{\partial t}\\
                                      &=& \rho_{f} \left(-\alpha_{f}  \dfrac{\partial T_{f}}{\partial t}+\beta_{f} \dfrac{\partial p}{\partial t} \right),                                 
\end{eqnarray}
 where $\alpha_{f}$ and $\beta_{f}$ are respectively the  thermal fluid expansivity
and its compressibility defined by
\begin{eqnarray}
 \alpha_{f}= -\dfrac{1}{\rho_{f}}\dfrac{\partial \rho_{f}}{\partial T_{f}},\,
\;\;\;\; \beta_{f}= \dfrac{1}{\rho_{f}}\dfrac{\partial \rho_{f}}{\partial p}.
\end{eqnarray}

Inserting equation \eqref{masscon3} and  \eqref{rockcomp}  in  equation \eqref{masscon2} yields
\begin{eqnarray}
\label{masscon4}
 -\phi \rho_{f}\alpha_{f}  \dfrac{\partial T_{f}}{\partial t}+ \rho_{f} \left( \phi \beta_{f} +  \phi_{0} \alpha_{b}\right)\dfrac{\partial p}{\partial t}&=&\nabla \cdot \left( \dfrac{\rho_{f}\mathbf{K}}{\mu} \left(\nabla p-\rho_{f}\mathbf{g}\right)\right)+Q_{f},
\end{eqnarray}
The state functions  $\mu$,$\rho_{f},\alpha_{f},\beta_{f}$ and $c_{pf}$  can  found \cite{finem,steamt}.

 The model problem is therefore to find the functions $(T_{s},T_{f},p)$ satisfying the nonlinear equations (\ref{heat}) and (\ref{masscon4}) subject to (\ref{darcy}). 
Notice if  $he\rightarrow \infty $ we have the equilibrium state of the heat with $T_{s}=T_{f}$.

\section{Numerical Schemes}
\label{sec3}
\subsection{Space discretization}
We use the finite volume method\cite{FV,FV1} on a structured mesh $\mathcal{T}$, 
with maximum mesh size $h$.
We denote by ($\Omega_{i}$) the family of control volumes of mesh $\mathcal{T}$.  The finite volume method space discretization consists in:
 \begin{enumerate}
  \item Integrate each equations of  \eqref{heat} and  \eqref{masscon4} over  each control volume $\Omega_{i}$.
  \item Use the divergence theorem to convert the volume integral into the surface integral in all divergence terms.
   \item Use two-point  flux approximations
 for diffusion  heat and flow fluxes \cite{FV}
\begin{eqnarray}
f_{i}^{1}=\int_{\partial \Omega_{i}}\mathbf{n}_{i} \cdot \left( \mathbf{k_{s}} \nabla T_{s}\right)ds,\;\;
f_{i}^{2}=\int_{\partial \Omega_{i}}\mathbf{n}_{i} \cdot \left(\mathbf{ k_{f}} \nabla T_{f}\right)ds,\;\;\\
f_{i}^{3}=\int_{\partial \Omega_{i}}\mathbf{n}_{i} \cdot \left(\dfrac{ \rho_{f}\mathbf{K}}{\mu}\nabla  p\right)ds.
\end{eqnarray}
\item Use the standard upwind weighting \cite{FV} for the convective (advective) flux
\begin{eqnarray}
f_{i}^{4}=\int_{\partial \Omega_{i}}\mathbf{n}_{i} \cdot \left( \rho_{f} c_{pf}T_{f} \mathbf{v}\right)ds.\;\;
\end{eqnarray}
\end{enumerate}
Here we denote by $\mathbf{n}_{i}$ the unit  normal vector  to $\partial \Omega_{i}$ outward  to $i$ 
 and  $ds$ the elementary surface measure. 
For an edge $\sigma$ of the control volume  $i$,
$\mathbf{n}_{i,\sigma}$  will denote  the  unit normal vector to $\sigma$ outward to $\Omega_i$.

Let us illustrate the spatial discretization of the second equation of \eqref{heat} on a structured mesh $\mathcal{T}$
(the two-point flux approximation is sufficient for so-called {$\bf K$}-orthogonal grids).
 We denote by $\;\mathcal{E}_{int} $ 
the set of interior  edges of the control volumes of $ \mathcal{T}$. For any function $X$,   $X_{i}(t)$  denotes
the approximation of $X$ at time $t$ at the center of the control 
volume $\Omega_i\in \mathcal{T}$ and $X_{\sigma}(t)$ the approximation of $X$ at time $t$ at 
the center of the edge $\sigma$.
 For a control volume $\Omega_i \in \mathcal{T}$,  we denote by $\mathcal{E}_{i}$ the set
 of edges of $\Omega_i$, 
$\vert \Omega_i \vert $ the Lebesgue measure of the control
 volume $\Omega_i\in \mathcal{T}$, $i\mid j$  the edge interface between the control volume $\Omega_i$ 
and the control volume $ \Omega_j\neq \Omega_i$, $d_{i,j}$  
the distance between the center of the control volume $\Omega_i$ and center of the control volume $\Omega_j$, 
and $d_{i,\sigma}$  the distance   between the control volume $\Omega_i$ and the edge  $\sigma$.
Letting $\mathbf{k_{s}}:= \mathbf{k}$, we  therefore have 
\begin{eqnarray*}
 f_{i}^{2}&=&\int_{\partial \Omega_{i}}\mathbf{n}_{i} \cdot \left( \mathbf{k} \nabla T_{f}\right)ds = \underset{\sigma \in\mathcal{E}_{i}} {\sum }f_{i,\sigma}\\
f_{i,\sigma}(t) &\approx&  \vert \sigma\vert\;k_{i,\sigma}\dfrac{T_{f\sigma}(t)-T_{fi}(t)}{d_{i,\sigma}},\,\,\,\,\,\,\,\,\sigma \in \partial \Omega \ \\
f_{i,\sigma}(t)&\approx&-\tau_{\sigma}\left( T_{fj}(t)-T_{fi}(t)\right), \,\,\,\,\,\,\,\,\sigma=i\mid j \in \mathcal{E}_{int} \\
\tau_{\sigma}&=&\vert \sigma\vert \dfrac{k_{i,\sigma}k_{j,\sigma}}{k_{i,\sigma}\; d_{i,\sigma}+k_{j,\sigma}
 d_{j,\sigma}} \\
k_{i,\sigma}&=&\vert \mathbf{k}_{i}\,\mathbf{n}_{i,\sigma}\vert,\quad\mathbf{k}_{i}=\dfrac{1}{\vert \Omega_i \vert }\int_{i}\mathbf{k}(\mathbf{x})d\mathbf{x},\\
\end{eqnarray*}
These approximations are for interior edges and Dirichlet Boundary condition. For a Neumann boundary 
condition, $\mathbf{n}_{i} \cdot \left( \mathbf{k} \nabla T_{f}\right)$ is naturally given.
In the case of a discrete-fracture model, we make adjustments to the spatial discretization following the approach in \cite{torh,TPFAdfm}.

To determine the convective fluxes, we set $X=\rho_{f} c_{pf}T_{f} $. Standard upwind weighting  yields
\begin{eqnarray}
f_{i}^{4}&=&\int_{\partial \Omega_{i}}\mathbf{n}_{i} \cdot \left( \rho_{f} c_{pf}T_{f} \mathbf{v}\right) ds = \underset{\sigma \in\mathcal{E}_{i}} {\sum }\int_{\sigma} \mathbf{n}_{i,\sigma} \cdot X \mathbf{v} ds
\approx \underset{\sigma \in\mathcal{E}_{i}}{\sum }  v_{i,\sigma} X_{\sigma,+}(t)\nonumber\\ 
 v_{i,\sigma}&=& \int_{\sigma}\mathbf{n}_{i,\sigma}\mathbf{v} ds\\
X_{\sigma,+}(t)&=&\left\lbrace \begin{array}{l}
  X_{i}( t)\quad \text{if} \quad v_{i,\sigma}\geqslant 0  \\  
       \newline\\             
   X_{j}( t) \quad  \text{if} \quad  v_{i,\sigma}< 0         
\end{array}\right.
\;\; \;\;\text{for}\;\;\;\sigma = i\mid j\\ 
X_{\sigma,+}(t)&=&\left\lbrace \begin{array}{l}
    X_{i}( t)\quad \text{if} \quad  v_{i,\sigma}\geqslant 0 \\
\newline\\
    X_{\sigma}(t) \quad  \text{if}  \quad  v_{i,\sigma}< 0                    
                             \end{array}\right.
\;\; \text{for}\;\;\sigma \in\mathcal{E}_{i}\cap \partial \Omega .
\end{eqnarray}

Reorganizing these space approximations yield the following system of ODEs
\begin{eqnarray}
\label{spaceheat}
\left\lbrace \begin{array}{rcl}
 \dfrac{ dT_{s}^{h}}{ dt}& = &G_{1}(T_{s}^{h},T_{f}^{h},t)\\
\newline\\
\dfrac{ dT_{s}^{h}}{ dt}&=&G_{2}(T_{s}^{h},T_{f}^{h},p_{h},t)\\
\newline\\
\dfrac{ dp_{h}}{ dt}&=&G_{3}( p_{h},T_{f}^{h},t)+ \dfrac{(\phi\alpha_{f})(T_{f}^{h},p_{h})}{(\phi\beta_{f}+\phi_{0} \alpha_{b})(T_{f}^{h},p_{h})} \cdot\dfrac{ dT_{s}^{h}}{ dt}
\end{array}\right.
\\
\Leftrightarrow 
\label{spaceheat1}
\left\lbrace \begin{array}{l}
\dfrac{ dT_{h}}{ dt}=G(T^{h},p_{h},t)\\
\newline\\
\dfrac{ dp_{h}}{ dt}=G_{3}( p_{h},T_{f}^{h},t)+\dfrac{(\phi\alpha_{f})(T_{f}^{h},p_{h})}{(\phi\beta_{f}+\phi_{0} \alpha_{b})(T_{f}^{h},p_{h})}\cdot G_{2}(T_{s}^{h},T_{f}^{h},p_{h},t)\\
\newline\\
 G(T^{h},p_{h},t)=(G_{1}(T_{s}^{h},T_{f}^{h},t),G_{2}(T_{s}^{h},T_{f}^{h}, p_{h}, t))^{T}\\
\newline\\
T_{h}=(T_{s}^{h},T_{f}^{h})^{T}\approx (T_{s},T_{f})^{T}.
\end{array}\right.
\end{eqnarray}
For a given initial pressure $p_{h}(0)$, with corresponding initial velocity $\mathbf{v}_{h}(0)$) and initial temperature $T_{h}(0)$, the technique used in the paper
consists in solving successively the systems
 \begin{eqnarray}
\label{spaceheat2}
\left\lbrace \begin{array}{l}
\dfrac{ dT_{h}}{ dt}=G(T^{h},p_{h},t)\\
\newline\\
T_{h}(0),\, p_{h}(0) \;\;\; \textbf{given}
\end{array}\right.
\end{eqnarray}
and 
\begin{eqnarray}
\label{spaceheat3}
\left\lbrace \begin{array}{l}
\dfrac{ dp_{h}}{ dt}=G_{3}( p_{h},T_{f}^{h},t)+   \dfrac{(\phi\alpha_{f})(T_{f}^{h},p_{h})}{(\phi\beta_{f}+\phi_{0} \alpha_{b})(T_{f}^{h},p_{h})}\cdot G_{2}(T_{s}^{h},T_{f}^{h},p_{h},t)=G_{4}(T_{h}^{h},p_{h},t) \\
\newline\\
T_{h}(0), \,p_{h}(0) \;\;\; \textbf{given}.
\end{array}\right.
\end{eqnarray}
The ODEs \eqref{spaceheat1} is usualy stiff since the smaller absolute value of 
the eigenvalues of the Jacobian matrix is usually closed to zero.
\subsection{Time discretizations}
We now present our numerical methods  for the ODEs \eqref{spaceheat1}
 based on the sequential approach.
Considering a sequential solution approach, we start by presenting low order time discretization 
schemes and their  stability properties before we give some higher order schemes.
\subsubsection{$\theta$-Euler Schemes}
We briefly describe the standard integrators 
that will be used for comparison with Exponential Rosenbrock-Euler method and Rosenbrock-type methods.
Consider the ODEs \eqref{spaceheat2} and \eqref{spaceheat3} within the 
 interval $[0, \tau],\,\tau>0$.  
Given a time-step $\tau_{n}$, applying the $\theta-$Euler scheme with respect to 
the function $T_{h}$ in the ODEs \eqref{spaceheat2} yields
\begin{eqnarray}
\label{temscheme}
 \left\lbrace \begin{array}{l}
   \dfrac{T_{h}^{n+1}-T_{h}^{n}}{\tau_{n}} = \theta G(T_{h}^{n+1},p_{h}^{n},t_{n+1})+(1-\theta)G(T_{h}^{n},p_{h}^{n},t_{n})\\
  \newline\\
    T_{h}(0)=T_{h}^{0},\;\;\;\;\; 0\leq \theta\leq 1.
\end{array}\right.
\end{eqnarray}
For $\theta \neq 0$ the scheme is implicit and 
given the approximate solutions $T_{h}^{n}$  and  $p_{h}^{n}$ at time $t_{n}$, the solution $T_{h}^{n+1}$ at time  $t_{n+1}$ is obtained by solving the nonlinear equation
\begin{eqnarray}
\mathcal{F}(X)= X- \tau_{n} \theta G(X,p_{h}^{n}, t_{n+1})- \tau_{n}(1-\theta)G(T_{h}^{n},p_{h}^{n},t_{n})-T_{h}^{n}=0,\,\,\,X=T_{h}^{n+1},
\end{eqnarray}
which is solved using the Newton method. For
efficiency, all linear systems are solved using the
Matlab function bicgstab with ILU(0) preconditioners with no fill-in, which are updated at each time step.

To solve the ODEs \eqref{spaceheat3} we apply again the $\theta$-Euler method, but with respect 
to $p_{h}$, yielding
\begin{eqnarray}
\label{prescheme}
 \left\lbrace \begin{array}{l}
   \dfrac{p_{h}^{n+1}-p_{h}^{n}}{\tau_{n}} = \theta G_{4}(T_{h}^{n+1},p_{h}^{n},t_{n+1})+(1-\theta)G_{4}(T_{h}^{n},p_{h}^{n},t_{n})\\
  \newline\\
    p_{h}(0)=p_{h}^{0}.
\end{array}\right.
\end{eqnarray}
For $\theta=1/2$ the  $\theta$-Euler scheme is second order in time, and for $\theta \neq \frac{1}{2}$ 
the scheme is first order in time.
In this paper, the standard sequential approach to solve the ODEs \eqref{spaceheat} consists of applying successively the schemes \eqref{temscheme}
and \eqref{prescheme}. 
 
\subsubsection{Exponential Rosenbrock-Euler Method and Rosenbrock-type Methods} 
To introduce the Exponential Euler-Rosenbrock Method (also called Exponential Euler Method), let us first consider the following system of ODEs
\begin{eqnarray}
\label{extraodes}
\left\lbrace \begin{array}{l}
 \dfrac{d\mathbf{y}}{dt}= \mathbf{L}\mathbf{y}+\mathbf{g}(\mathbf{y},t),\;\;\;\;\ t\in [0,\tau]\\
\mathbf{y}(0)=\mathbf{y}_{0},
\end{array}\right.
\end{eqnarray}
which appear after spatial discretisation of semilinear parabolic PDEs.
Here $\mathbf{L}$ is a stiff matrix and $\mathbf{g}$ a nonlinear function.
This allows us to write the exact solution
of \eqref{extraodes} as
\begin{eqnarray*}
\mathbf{y}(t_{n})=
e^{t_{n}\mathbf{L}}\mathbf{y}_{0}+\int_{0}^{
  t_{n}} e^{(t_{n}-s)\mathbf{L} }\mathbf{g}(\mathbf{y}(s),s)ds.
\end{eqnarray*}
Given the exact solution at the
time $t_{n}$, we can construct the corresponding solution at $t_{n+1}$
as
\begin{eqnarray}
\label{exactnu}
\mathbf{y}(t_{n+1})&=&e^{\tau_{n}\mathbf{L}}\mathbf{y}(t_{n})+\int_{0}^{ \tau_{n}} e^{(\tau_{n}-s)\mathbf{L} }\mathbf{g}(\mathbf{y}(t_{n}+s),s+t_{n})ds,\\
\tau_{n}&=& t_{n+1}-t_{n}.
\end{eqnarray}
Note that the expression in \eqref{exactnu} is still an exact
solution. The idea behind Exponential Time Differencing (ETD) is to
approximate $\mathbf{g}(\mathbf{y}(t_{n}+s),t_{n}+s)$
by a suitable polynomial~\cite{MC}. The simplest case is 
when $\mathbf{g}(\mathbf{y}(t_{n}+s),t_{n}+s)$ is
approximated by the constant $
\mathbf{g}(\mathbf{y}(t_{n}),t_{n})$. The corresponding (ETD1) scheme is given by
\begin{equation}
\label{etd}
\mathbf{y}^{n+1}  = e^{\tau \mathbf{L}}
\mathbf{y}^{n}+\tau_{n}\varphi_{1} (\tau_{n} \mathbf{L})\mathbf{g}(
\mathbf{y}^{n},t_{n}),
\end{equation}
where
$$
\varphi_{1}(A)= \underset{i=1}{\sum^{\infty}} \dfrac{A^{i-1}}{i!}= A^{-1}(e^{A}-I)\;\;( \text{if $A$ is invertible}).
$$
Note that the ETD1 scheme in  \eqref{etd} can be rewritten as
\begin{eqnarray}
\label{etd1}
\mathbf{y}^{n+1}  =  \mathbf{y}^{n}+\tau_{n}\varphi_{1} (\tau_{n} \mathbf{ L})(\mathbf{ L}\mathbf{y}^{n}+\mathbf{g}(\mathbf{y}^{n},t_{n})).
\end{eqnarray}

This new expression has the advantage that it is  computationally more
efficient as only one matrix exponential function needs to be
evaluated at each step.

Recently, the ETD1 scheme was applied to advection-dominated reactive transport in heterogeneous
 porous media  \cite{ATthesis,Antoine}. For this problem, a rigorous convergence proof is established for the case of a finite volume discretization in space 
 \cite{ATthesis,tambueb}. In these works, it was observed that the
exponential methods were generally more accurate and efficient than standard implicit methods. 

As our systems \eqref{spaceheat2} and \eqref{spaceheat3} are nonlinear, we need to linearize before applying the 
ETD1 scheme. 
Consider the system the ODEs \eqref{spaceheat2}. For simplicity we assume that it is autonomous
\begin{eqnarray}
 \left\lbrace \begin{array}{l}
\dfrac{d T_{h}}{ dt} = G(T_{h},p_{h})\;\;\;\;\;\;\;
\newline\\
T_{h}(0)=T_{0}.
\end{array}\right.
\end{eqnarray}

Let $T_{h}^{n}$ and $p_{h}^{n}$ be the numerical approximations of the exact solutions $T(t_{n})$ and  $p(t_{n})$. To obtain the numerical approximation
$T_{h}^{n+1}$  of the exact solution $T(t_{n+1})$, we linearize  $G(T_{h},p_{h}^{n})$ at $T_{h}^{n}$  and obtain the following semilinear ODEs
\begin{eqnarray}
\label{linear}
 \dfrac{d T_{h}}{dt}&=&\mathbf{J}_{n}T_{h}(t)+\mathbf{g}(T_{h}(t),p_{h}^{n})\;\;\; t_{n}\leq t\leq t_{n+1},
\end{eqnarray}
where $\mathbf{J}_{n}$ denotes the Jacobian of the function $G$ respect to $T_{h}$ and $\mathbf{g}$ the remainder given by
\begin{eqnarray}
 \mathbf{J}_{n}=D_{T_{h}}G(T_{h}^{n},p_{h}^{n}),\;\;\;\;\mathbf{g}(T_{h}(t),p_{h}^{n})=G(T_{h}(t),p_{h}^{n})-\mathbf{J}_{n}T_{h}(t).
\end{eqnarray}
Applying the ETD1 scheme to \eqref{linear} yields
\begin{eqnarray}
\label{EERM}
 T_{h}^{n+1}=T_{h}^{n}+ \tau_{n}\varphi_{1}(\tau_{n}\mathbf{J}_{n})G(T_{h}^{n},p_{h}^{n}),\;\;\; \tau_{n}=t_{n+1}-t_{n}.
\end{eqnarray}
This scheme, called the Exponential Rosenbrock-Euler method (EREM)\cite{calarioster} (or  the Exponential Euler method (EEM))  has been
 reinvented in  different names (see references in \cite{carr}). 

The EREM scheme is second order in time  \cite{calarioster,carr} for autonomous problems.
To deal with  non-autonomous problems, while conserving the second order accuracy of
 EERM scheme, it must first be converted to autonomous problems. The corresponding version is given in the next section by equation\eqref{nEERMg}.

The scheme \eqref{EERM} contains the exponential matrix  function $\varphi_{1}$. To obtain the simplest Rosenbrock-type methods, the exponential function 
is approximated by the  following rational function
\begin{eqnarray}
\label{ratapp}
\varphi_{1} (\tau_{n}\mathbf{J}_{n}) \approx \left( \mathbf{I}-\tau_{n} \gamma \mathbf{J}_{n}\right)^{-1}, 
\end{eqnarray}
where $\gamma> 0$  is a parameter.
 For a given parameter $\gamma> 0$, using the approximation \eqref{ratapp}  in equation \eqref{EERM}, 
the corresponding  Rosenbrock-type Method (also called linear implicit methods)
 is given  by
\begin{eqnarray}
\label{ROS}
 T_{h}^{n+1}=T_{h}^{n}+ \tau_{n}\left( \mathbf{I}-\tau_{n} \gamma \mathbf{J}_{n}\right)^{-1} G(T_{h}^{n},p_{h}^{n}).\;
\end{eqnarray}


For the parameter $\gamma=1/2$, the corresponding Rosenbrock-type Method is 
 order two in time for regular solutions and order 1 if  $\gamma\neq 1/2$.

To  solve the ODEs \eqref{spaceheat3} we apply again the EREM scheme or  the Rosenbrock-type method \eqref{ROS}, but with respect 
to $p_{h}$, which yields respectively
\begin{eqnarray}
\label{EERM1}
\left\lbrace \begin{array}{l}
 p_{h}^{n+1}=p_{h}^{n}+ \tau_{n}\varphi_{1}(\tau_{n}\mathcal{J}_{n})G_{4}(T_{h}^{n},p_{h}^{n}),\\
\newline\\
\mathcal{J}_{n}=D_{p_{h}}G_{4}(T_{h}^{n},p_{h}^{n}),
\end{array}\right.
\end{eqnarray}
and
\begin{eqnarray}
\label{ROS1}
 p_{h}^{n+1}=p_{h}^{n}+ \tau_{n} \left( \mathbf{I}-\tau_{n} \gamma \mathcal{J}_{n}\right)^{-1} G_{4}(T_{h}^{n},p_{h}^{n}).
\end{eqnarray}
As with $\theta-$Euler methods, the sequential approaches with the Exponential Rosenbrock-Euler method and  Rosenbrock-type methods 
proposed in this paper consists in solving the ODEs \eqref{spaceheat}  by 
applying successively the schemes \eqref{EERM} and  \eqref{ROS} to the  ODEs \eqref{spaceheat2} 
and  the schemes \eqref{EERM1} and  \eqref{ROS1} to the  ODEs \eqref{spaceheat3}.

The sequential technique presented is the so called Trotter splitting and is generally first order accurate\cite[p. 7]{lie}.
The alternative technique is the so called Strang splitting which consists 
to apply the schemes \eqref{EERM} and  \eqref{ROS} to the  ODEs \eqref{spaceheat2} with time step $\dfrac{\tau_{n}}{2}$,
 afterward apply the schemes \eqref{EERM1} and  \eqref{ROS1} to the ODEs \eqref{spaceheat3} with time step $\tau_{n}$ and 
then apply again the schemes \eqref{EERM} and  \eqref{ROS} to the  ODEs \eqref{spaceheat2} with time step $\dfrac{\tau_{n}}{2}$. Strang splitting is 
formally second -order accurate in time for sufficiently smooth solution \cite[p. 7]{lie}. 
 We can obviously observe  that  this approach is less efficient than Trotter splitting.
Trotter splitting and Strang splitting are called multiplicative operator splittings. 

In the sequel, sequential approach will mean Trotter splitting.

\subsubsection{Stability properties of numerical schemes and  higher order Rosenbrock-type methods}
One of the important features of any numerical scheme is its stability properties.
Our goal here is to study the stability properties of the schemes presented in  the previous 
section and high order Rosenbrock-type methods. A special interest will be given to two Rosenbrock-type methods 
of order two and three because of their good stability properties. In applying the high order Rosenbrock-type methods, we will
use the previously presented sequential approaches to solve the ODEs \eqref{spaceheat2} and \eqref{spaceheat3}.

Consider the following ODEs
\begin{eqnarray}
\label{extraodes1}
\left\lbrace \begin{array}{l}
 \dfrac{d\mathbf{y}}{dt}= \mathbf{f}(\mathbf{y},t),\;\;\;\;\ t\in [0,\tau]\\
\mathbf{y}(0)=\mathbf{y}_{0},
\end{array}\right.
\end{eqnarray}
where $\mathbf{f}$ is nonlinear function.   The corresponding $\theta-$Euler scheme is given by
\begin{eqnarray}
 \left\lbrace \begin{array}{l}
   \dfrac{\mathbf{y}^{n+1}-\mathbf{y}^{n}}{\tau_{n}} = \theta  \mathbf{f}(\mathbf{y}^{n+1},t_{n+1})+(1-\theta)\mathbf{f}(\mathbf{y}^{n},t_{n})\\
  \newline\\
   \mathbf{y}(0)=\mathbf{y}_{0}\;\;\;\;\; 0\leq \theta\leq 1.
\end{array}\right.
\end{eqnarray}
Note that the exponential Rosenbrock--Euler method presented in the previous section is for autonomous ODEs.
 Before applying it to non-autonomous system \eqref{extraodes1}, transformation  $\hat{\mathbf{y}} =\left(\mathbf{y},t\right)^{T}$ 
must be performed to obtain autonomous ODEs. Given the numerical solution
$\hat{\mathbf{y}}^{n}=(\mathbf{y}^{n},t_{n})^{T}$, the linearization equation leading to $\hat{\mathbf{y}}^{n+1}$ is given by
\begin{eqnarray*}
\hat{\mathbf{y}}'(t)=
\left[ \begin{array}{c}
        \mathbf{y}'(t)\\
          t'
       \end{array}
\right]
&=&
\left[ \begin{array}{cc}
        D_{\mathbf{y}}\mathbf{f}(\mathbf{y}^{n},t_{n}) & D_{t}\mathbf{f}(\mathbf{y}^{n},t_{n})\\
          0 & 0\\
       \end{array}
\right]
\left[ \begin{array}{c}
        \mathbf{y}(t)\\
          t
       \end{array}
\right]
+
\left[ \begin{array}{c}
       \mathbf{g}_{n} (\mathbf{y}(t),t)\\
          1
       \end{array}
\right]
\newline\\
\mathbf{g}_{n} (\mathbf{y}(t),t)&=&\mathbf{f}(\mathbf{y}(t),t)-D_{\mathbf{y}}\mathbf{f}(\mathbf{y}^{n},t_{n})\mathbf{y}(t)-D_{t}\mathbf{f}\, t.
\end{eqnarray*}
 Using this transformation and \cite[Lemma 1]{calarioster}, the corresponding EREM scheme  for non-autonomous system is given by  
\begin{eqnarray}
\label{nEERMg}
\left\lbrace \begin{array}{l}
 \mathbf{y}^{n+1}=\mathbf{y}^{n}+ \tau_{n}\varphi_{1}(\tau_{n}\mathcal{J}_{n})\mathbf{f}(\mathbf{y}^{n},t_{n})+ \tau_{n}^{2}\varphi_{2}(\tau_{n}\mathcal{J}_{n})D_{t}\mathbf{f}(\mathbf{y}^{n},t_{n})\\
\newline\\
\mathcal{J}_{n}=D_{\mathbf{y}}\mathbf{f}(\mathbf{y}^{n},t_{n}).
\end{array}\right.
\end{eqnarray}

The exponential functions $\varphi_{i}$ are defined by 
\begin{eqnarray}
 \varphi_{i}(z)= \int_{0}^{1} e^{(1-s)z}\dfrac{s^{i-1}}{(i-1)!}ds, \qquad  \qquad   \qquad  \,\,\,\, i \geq 1.
\end{eqnarray}
These functions satisfy the recurrence relations
\begin{eqnarray}
 \varphi_{i}(z)=\dfrac{\varphi_{i-1}(z)-\varphi_{i-1}(0)}{z},\;\,\;\;\,\,\,\,\varphi_{0}=e^z.
\end{eqnarray}

The corresponding lower order  Rosenbrock-type methods is given  by
\begin{eqnarray}
\label{ROSf}
 \mathbf{y}^{n+1}=\mathbf{y}^{n}+ \tau_{n} \left( \mathbf{I}-\tau_{n} \gamma \mathcal{J}_{n}\right)^{-1}\left[\mathbf{f}(\mathbf{y}^{n},t_{n})+\gamma \tau_{n} D_{t} \mathbf{f}(\mathbf{y}^{n},t_{n})\right].
\end{eqnarray}
  
In order to study their stability properties, we apply the $\theta-$Euler, EREM  and ROSM  schemes to the linear ODEs
 $\mathbf{y}' = \lambda \mathbf{y}$
  with  constant time step $\Delta t$. We therefore have 
\begin{eqnarray}
 \mathbf{y}^{n+1}=R(z)\mathbf{y}^{n},\,\,\,\  z= \Delta t\lambda\,
\end{eqnarray}
where  
 $ R(z) =\dfrac{1+z(1-\theta)}{1-\theta z} $ for  the $\theta-$Euler scheme, $R(z)=\exp{(z)}$ for the EREM  
scheme  and $ R(z) =\dfrac{1+z(1-\gamma)}{1- \gamma z} $ for the ROSM scheme.

The function $R(z)$ is called the stability function of the method. The set
$$S=\left\lbrace z \in \mathbb{C},\,\,\,\ \vert R(z)\vert \leq 1 \right\rbrace $$

is called the stability domain of the method.  A numerical method is  A -stable if its stability domain $S$ satisfies
$$ S\supset \mathbb{C}^{-}=\left\lbrace z\in \mathbb{C},\,\,\,\ \text {Re}\, z \leq 0 \right\rbrace.$$ 
Let us study the A-stability of the $\theta-$Euler. Let $z=x+iy \in \mathbb{C}^{-}$, we have
\begin{eqnarray*}
 \vert R(z)\vert \leq 1 & \Leftrightarrow &\vert 1+(1-\theta )z\vert \leq \vert 1-\theta z\vert\\
& \Leftrightarrow & (1+(1-\theta)x)^{2} +((1-\theta)y)^{2} \leq ((1- \theta)x)^{2}+(\theta y)^{2}\\
& \Leftrightarrow & (1-2\theta)(x^{2}+y^{2})+2x \leq 0.
\end{eqnarray*}

We can therefore observe that the $\theta-$Euler
 scheme is A -stable if $$1-2\theta \leq 0 \Leftrightarrow \theta \geq 1/2,$$
ROSM scheme is  A -stable if $\gamma \geq 1/2$ and EREM scheme is  A -stable.

A-stability is not the whole answer to the problem of stiff equations, excellent numerical methods for super-stiff equations would be L-stable.

Numerical methods are L-stable if they are  A-stable  and in addition (see \cite{hairer})
$$ \underset{z\rightarrow -\infty} {\lim}\vert R(z)\vert = 0.$$  We can  observe that the $\theta-$Euler scheme is L-stable if $\theta =1$, 
the ROSM scheme is L-stable if $\gamma=1$  and  the  EREM scheme is  L-stable.

In the sequel we will  use the ROSM schemes with $\gamma=1$ and $\gamma=1/2$, which will be denoted respectively by ROSM$(1)$ and  ROSM$(1/2)$.

The s-stage  Rosenbrock-type methods for the ODE \eqref{extraodes1} are given by 
\begin{eqnarray}
\label{ROSMh}
 \left\lbrace \begin{array}{l}
\left( \dfrac{1}{\tau_{n}  \gamma}\mathbf{I}- \mathcal{J}_{n}\right)\mathbf{k}_{n i}=\mathbf{f}(\mathbf{y}_{n}+ \underset{j=1}{\sum ^{i-1}}a_{i j}\mathbf{k}_{n j},t_{n}+\alpha_{i} \tau_{n})
   - \underset{j=1}{\sum ^{i-1}}\dfrac{c_{i j}}{\tau_{n}}\mathbf{k}_{n j}\\
 \qquad  \qquad \qquad \qquad+ \tau_{n}  \gamma_{i} D_{t} \mathbf{f}(\mathbf{y}^{n},t_{n})\,\,\,\,\,\,\, i=1, \cdots, s \\
\newline\\
    \qquad \qquad \quad  \mathbf{y}^{n+1}=\mathbf{y}^{n}+ \underset{i=1}{\sum ^{s}}b_{i}\mathbf{k}_{n i}\\
\newline\\
    \qquad \qquad \quad  \mathbf{y}_{1}^{n+1}=\mathbf{y}^{n}+ \underset{i=1}{\sum ^{s}}\hat{b}_{i}\mathbf{k}_{n i}.
\end{array}\right.
\end{eqnarray}
The  coefficients $a_{i j},c_{i j},\alpha_{i},\gamma, \gamma_{i},b_{i}$ 
are obtained by using  the consistency conditions required to achieve  the desirable  order  of convergence $p$ in time. 
Different ways to find these coefficients are presented in the literature (see \cite{Rosen,hairer,langb,lang} and references therein).
The approximation $\mathbf{y}_{1}^{n+1}$ is called an embedded approximation associated to 
the  s-stage  Rosenbrock-type approximation $\mathbf{y}^{n+1}$ and is used to control the local errors for adaptivity purpose. 
The coefficients $\hat{b}_{i}$ are determined  using the consistency conditions such that  the embedded approximation is order $p-1$.
In  this work,  we use  the second  order scheme ROS2(1), where the coefficients are given in \tabref{ros2}, 
and also the third order scheme denoted ROS3p in \cite{lang}, which uses  additional conditions to avoid order reduction (see \tabref{ros3}).
 
The ROS2(1) scheme is L-stable and the  ROS3p scheme is A-stable (see \cite{langb}).
\begin{table}[h!]
\begin{center}
\begin{tabular}{|l|l|}
\hline
 $\gamma=1.707106781186547e^{+00} $& \\
\hline
$a_{11}=0 $ & $c_{11}=5.857864376269050 e^{-01}$   \\
$a_{21}=5.857864376269050 e^{-01} $ & $c_{21}=1.171572875253810 e^{+00}$   \\
$a_{22}=0 $ &   $c_{22}= 5.857864376269050 e^{-01}$ \\
\hline
 $\gamma_{1}=1.707106781186547 e^{+00} $&  $\alpha_{1}=0$\\
$\gamma_{2}=-1.707106781186547 e^{+00} $&  $\alpha_{2}=1$\\
\hline
 $b_{1}=8.786796564403575e^{-01} $&   $ \hat{b}_{1}=5.857864376269050 e^{-01}  $ \\
$b_{2}=2.928932188134525e^{-01} $ &  $ \hat{b}_{2}=0$\\
\hline
\end{tabular}
\end{center}
\caption{Coefficients of the ROS2(1) scheme from \cite{langb}.}
\label{ros2}
\end{table}

\begin{table}[h!]
\begin{center}
\begin{tabular}{|l|l|}
\hline
 $\gamma=7.886751345948129e^{-01} $& \\
\hline
$a_{11}=0 $ & $c_{11}=1.267949192431123 e^{+00}$   \\
$a_{21}=1.267949192431123 e^{+00} $ & $c_{21}=1.607695154586736 e^{+00}$   \\
$a_{22}=0 $ &   $c_{22}= 1.267949192431123e^{+00}$ \\
$a_{31}= 1.267949192431123e^{+00}$ &   $c_{31}= 3.464101615137755e^{+00}$ \\
$a_{32}= 0$ &   $c_{32}= 1.732050807568877e^{+00}$ \\
$a_{33}= 0$ &   $c_{33}= 1.267949192431123e^{+00}$ \\
\hline
 $\gamma_{1}=7.886751345948129 e^{-01} $&  $\alpha_{1}=0$\\
$\gamma_{2}=-2.113248654051871 e^{-01}$&  $\alpha_{2}=1$\\
$\gamma_{3}=-1.077350269189626e^{+00}$&  $\alpha_{3}=1$\\
\hline
 $b_{1}=2$&   $ \hat{b}_{1}=2.113248654051871 e^{+00} $ \\
$b_{2}=5.773502691896258e^{-01} $ &  $ \hat{b}_{2}=1$\\
$b_{3}=4.226497308103742e^{-01} $ &  $ \hat{b}_{3}=4.226497308103742e^{-01}$\\
\hline
\end{tabular}
\end{center}
\caption{Coefficients of the ROS3p scheme from \cite{langb,lang}.}
\label{ros3}
\end{table}
The implementation of Rosenbrock-type schemes is straightforward as there are no nonlinear equations to solve at each time step. For
efficiency, all linear systems are solved using the
Matlab function bicgstab with ILU(0) preconditioners with no fill-in. 
The time step adaptivity can be performed using  the standard error control  and the step size prediction
as in \cite[pp.112]{hairer} with an appropriate norm of $ \mathbf{y}^{n+1}-\mathbf{y}_{1}^{n+1}$.

\section{Implementation of Exponential Rosenbrock-Euler Method}
The key element in all exponential integrator schemes is the
computation of the matrix exponential functions, the so called  $\varphi_{i}-$
functions. There are many techniques available for that task\cite{osterreview}.
Standard Pad\'{e} approximation compute
at every time step the whole matrix exponential functions and are therefore memory and time consuming for large problems. Krylov subspace technique
  and the real fast Leja points technique are proved 
to be  efficient for this computation for large systems \cite{LE,carr,TLGspe1,kry,SID,ATthesis}. Let us summarize these techniques 
while solving ODE \eqref{extraodes1} with the EREM scheme.

\subsection{Krylov subspace technique}
The main idea of the Krylov subspace technique is to approximate the
action of the exponential matrix function  $\varphi_{i}( \tau_{n}\mathcal{J}_{n})$  on a
vector $\mathbf{v}$ by projection onto a small Krylov subspace
$K_{m} =\text{span} \left\lbrace
 \mathbf{v} ,\mathcal{J}_{n} \mathbf{v},\ldots,  \mathcal{J}_{n}^{m-1}
  \mathbf{v}\right\rbrace $~\cite{SID}. The approximation is formed
using an orthonormal  basis of $\mathbf{V}_{m} =
\left[\mathbf{v}_{1} , \mathbf{v}_{2},\ldots, \mathbf{v}_{m}
\right] $ of the Krylov subspace $K_{m}$ and of its completion
$\mathbf{V}_{m+1}=\left[\mathbf{V}_{m},\mathbf{v}_{m+1}\right]$. The
basis is found by Arnoldi iteration, 
which uses stabilized Gram-Schmidt to produce a
sequence of vectors that span the Krylov subspace
(see \algref{alg:al2}). 

Let $\mathbf{e}_{i}^{j}$ be the $i^{\textrm{th}}$ standard basis vector of $\mathbb{R}^{j}$.
 We approximate $\varphi_{i}(\tau_{n} \,\mathcal{J}_{n})\mathbf{v}$ by
\begin{eqnarray}
\label{approphi}
  \varphi_{i}(\tau_{n}\,\mathcal{J}_{n})\mathbf{v}
  &\approx & \Vert \mathbf{v} \Vert_{2}\mathbf{V}_{m+1}\varphi_{i} ( \tau_{n}
  \, \overline{\mathbf{H}}_{m+1})\mathbf{e}_{1}^{m+1}
  \label{app}
\end{eqnarray}
with
$$
\overline{\mathbf{H}}_{m+1} = \left( \begin{array}{cc}
    \mathbf{H}_{m}&\b0\\
    0,\cdots,0, h_{m+1,m} & 0
  \end{array}\right) \qquad \text{where} \qquad
\mathbf{H}_{m} = \mathbf{V}_{m}^{T} \mathcal{J}_{n}\mathbf{V}_{m}=[h_{i,j}].
$$
The coefficient $h_{m+1,m}$ is recovered in the last iteration of
Arnoldi's iteration in \algref{alg:al2}. We denote by $\Vert.\Vert_{2}$ the standard Euclidan norm.
The approximation \eqref{approphi} is the first two terms of the  expansion given in \cite[Theorem 2]{SID}.
\begin{algorithm}
  \begin{algorithmic}[1]
 \STATE Initialise: $\mathbf{v}_{1} =\dfrac{\mathbf{v}}{\Vert \mathbf{v}\Vert_{2}}$ \COMMENT{normalisation}
\FOR{$j = 1 \cdots m$}
\STATE $\mathbf{w} = \mathcal{J}\mathbf{v}_{j}$
\FOR{$ i = 1 \cdots j $}
\STATE $h_{i,j} = \mathbf{w}^{T}\mathbf{v}_{i}$ \COMMENT{compute inner product to build  elements of the matrix $\mathbf{H}$}
\STATE $\mathbf{w} = \mathbf{w} - h_{i,j} \mathbf{v}_{i}$ \COMMENT{Gram--Schmidt process}
\ENDFOR
\STATE $ h_{j+1,j} = \Vert\mathbf {w}\Vert_{2} $
\STATE $\mathbf{v}_{j+1} = \dfrac{\mathbf{w}}{\Vert\mathbf{w}\Vert_{2} }$
\COMMENT{normalisation}
\ENDFOR
\end{algorithmic}
\caption{: Arnoldi's algorithm}
\label{alg:al2}
\end{algorithm}
For a small Krylov subspace (i.e, $m$ is small) a standard Pad\'{e}
approximation can be used to compute
$\varphi_{i}(\tau_{n}\overline{\mathbf{H}}_{m+1})$, but an efficient way 
used in ~\cite{SID} is to recover
$\varphi_{i}(\tau_{n}\overline{\mathbf{H}}_{m+1})\mathbf{e}_{1}^{m+1}$  
directly from the Pad\'{e} approximation of the exponential of a matrix
related to $\mathbf{H}_{m}$. 

Notice that this implementation can be done without explicit computation of the Jacobian matrix $\mathcal{J}_{n}$ as 
the Krylov subspace $K_{m}$  can be formed  by using  the  following approximations
\begin{eqnarray*}
\mathcal{J}_{n}\mathbf{v} &\approx& \dfrac{\mathbf{f}(\mathbf{y}^{n}+\epsilon \mathbf{v},t_{n})-\mathbf{f}(\mathbf{y}^{n},t_{n})}{\epsilon}\\
\,\,&\text{or}&\;\; \\
\mathcal{J}_{n}\mathbf{v}&\approx& \dfrac{\mathbf{f}(\mathbf{y}^{n}+\epsilon \mathbf{v},t_{n})-\mathbf{f}(\mathbf{y}^{n}-\epsilon \mathbf{v},t_{n})}{\epsilon},
\end{eqnarray*}
for a suitably chosen perturbation of $\epsilon$ (see \cite{JAF}), while solving the ODE \eqref{extraodes1}. These approximations
prove that  the Exponential Rosenbrock-Euler scheme with the Krylov subspace technique can be implemented using the free Jacobian technique.
The implementations in Expokit~\cite{SID} (for  function $\varphi_1$) and in \cite{P5} use the  truncation error in the approximation \eqref{approphi}  
to build the local error estimates (see \cite[Theorem 2]{SID}). 
The time step subdivisions depend on the given tolerance and the local errors.
\subsection{The real fast L\'{e}ja point technique}
This technique has been successfully applied to the  nonlinear advection diffusion-reaction-equation in \cite{LE,calarioster,LE2,parallelleja,ATthesis,Antoine} where advection plays a key role.  We will used it to solve the temperature equation \eqref{spaceheat2}.  
The key points of this method are as follows:
For a given vector $\mathbf{v}$, real fast \Leja points approximate $\varphi_{i}(\tau_{n} \,\mathcal{J}_{n})\mathbf{v} $  by $P_{m}(\tau_{n} \,\mathcal{J}_{n})\mathbf{v}$,
 where $P_{m}$ is an interpolation polynomial of  degree $m$ of $\varphi_{i}$
 at the sequence of points $\left\lbrace \xi_{i} \right\rbrace_{i=0}^{m} $ called spectral real fast \Leja points. These points $\left\lbrace \xi_{i} \right\rbrace_{i=0}^{m}$
 belong to the spectral focal interval  $\left[ \alpha, \beta\right]$  of the matrix  $\tau_{n}\mathcal{J}_{n}$, i.e. the focal interval of the
smaller ellipse containing all the eigenvalues of  $\tau_{n}\,\mathcal{J}_{n}$. This spectral interval can be estimated by  the well known  Gershgorin circle theorem~\cite{GC}. 
It has been shown that as the degree of the polynomial increases, and hence the number of
\Leja points increases, superlinear convergence
is achieved~\cite{LE2}; i.e.,
\begin{eqnarray}
  \underset {m \rightarrow \infty}{\lim}\Vert\varphi_{i}(\tau_n\mathcal{J}_{n})\mathbf{v}- P_{m}(\tau_n \,\mathcal{J}_{n})\mathbf{v}\Vert_{2}^{1/m}=0.
\end{eqnarray}
Set $\xi_{0}=\beta$, the sequence of
 fast L\'{e}ja points is generated by
\begin{eqnarray}
\label{lejaa}
  \underset {k=0}{\prod^{j-1}}\vert \xi_{j}-\xi_{k} \vert =\underset {\xi \in
 \left[\alpha, \beta \right] } {\max} \underset {k=0}{\prod^{j-1}}\mid \xi-\xi_{k} \mid \quad j=1,2,3,\cdots.
\end{eqnarray}
Given the Newton's form of the interpolating polynomial, $P_{m}$ is given by
\begin{eqnarray}
\label{approle}
 P_{m}(z)=\varphi_{i}\left[ \xi_{0}\right] + \underset{j=1}{\sum ^{m}}\varphi_{i}\left[ \xi_{0},\xi_{1}, \cdots,\xi_{j}\right] \underset{k=0}{\prod^{j-1}}\left( z-\xi_{k} \right)
\end{eqnarray}
where the divided differences $\varphi_i[\bullet]$ are defined recursively by
\begin{eqnarray}
\label{div}
\left\lbrace \begin{array}{l}
d_{0}=\varphi_{i}\left[ \xi_{0}\right] :=\varphi_{i}(\xi_{0}),\;\;
d_{1}=\varphi_{i}\left[ \xi_{0},\xi_{1}\right]:=\dfrac{\varphi_{i}\left[\xi_{1}\right]-\varphi_{i}\left[\xi _{0}\right]}{\xi_{1}-\xi_{0}}.\\
d_{i}= \varphi_{i}\left[ \xi_{0},\xi_{1},\cdots,\xi_{i}\right]=\dfrac{\varphi_{i}\left[\xi_{0},\xi_{1},\cdots,\xi_{i-2}, \xi_{i} \right]-\varphi_{i}\left[\xi_{0},\xi_{1},\cdots,\xi_{i-2},\xi_{i-1} \right]}{x_{i}-x_{i-1}}.
\end{array}\right.
\end{eqnarray}
Due to cancelation errors, this standard procedure cannot produce accurate divided differences
with magnitude smaller than machine precision.

It can be shown  \cite{Reichel} that the
divided differences of a function $f(c+\gamma \xi), \,c=(\alpha+\beta)/2, \,\gamma=(\beta-\alpha)/4 $ of the
independent variable $\xi$  at the  points 
$\left\lbrace \xi_{i}\right\rbrace_{i=0}^{m} \subset [-2, 2] $ are the
first column of the matrix function  $f(\mathbf{L}_{m})$, where
$$
\mathbf{L}_{m} =\tau_{n}\left(c \mathbf{I}_{m+1}+\gamma\mathbf{\widehat{L}}_{m}\right),\;\;\;
\mathbf{\widehat{L}}_{m}= \left( \begin{array}{cccccc}
    \xi_{0} & & & &&\\
    1& \xi_{1}& & &&\\
    & 1& \ddots&&&\\
    & & \ddots&\ddots&&\\
    & & &\ddots&\ddots\\
    & & &&1& \xi_{m}\\
  \end{array}\right). 
$$
 Here $\mathbf{I}_{m+1}$ is the identity matrix.
To compute $ (d_{i})_{0}^{m}=\varphi_{i}(\mathbf{L}_{m})e_{1}^{m+1}$, where  $e_{1}^{m+1}$ is
the first standard basis vector of $\mathbb{R}^{m+1}$, we apply Taylor
expansion of order $p$ with scaling and squaring or  Pad\'{e} approximation \cite{SID,ATthesis}.
The interest of the Newton interpolation comes from the fact that the approximation with a polynomial of degree $m$ is directly obtained
from the approximation with a polynomial of degree $m-1$, in fact  we have  
\begin{eqnarray}
\label{newt}
\left\lbrace \begin{array}{l}
  P_{m}(\tau_{n}\mathcal{J}_{n})\mathbf{v}=P_{m-1}(\tau_{n}\mathcal{J}_{n})\mathbf{v}+d_{m}q_{m}\\
q_{m}=\left((\mathcal{J}_{n}-c\mathbf{I})/\gamma -\xi_{m-1}\mathbf{I}\right)q_{m-1}\\
 P_{0}(\tau_{n}\mathcal{J}_{n})\mathbf{v}=d_{0}q_{0}\;\,\,\,\,\,q_{0}=\mathbf{v}.
\end{array} \right.
\end{eqnarray}
The error estimate from this approximation is given by
\begin{eqnarray*}
 \Vert e_{m}\Vert& =& \Vert P_{m+1}(\tau_n \mathcal{J}_{n})\mathbf{v}-P_{m}(\tau_{n}\mathcal{J}_{n})\mathbf{v}\Vert= \Vert d_{m}q_{m}\Vert\\
 &\approx& \Vert P_{m}(\tau_{n}\mathcal{J}_{n})\mathbf{v}-
\varphi_{i}(\tau_{n}\mathcal{J}_{n})\mathbf{v}\Vert,
\end{eqnarray*}
where $\Vert. \Vert $ is the weighted and scaled norm defined by 
\begin{eqnarray*}
 \Vert e_{m}\Vert =\sqrt{\dfrac{1}{N} \underset{i=1}{\sum^{N}} \left(\dfrac{e_{mi}}{\text{scal}_{i}} \right)^{2}}, \;\; \text{scal}_{i}= \text{tol}_a +  \text{tol}_r \cdot \Vert \mathbf{y}^{n} \Vert_{\infty},
\end{eqnarray*}
where $\text{tol}_a $ and  $\text{tol}_r $ denote respectively the desirable absolute and relative tolerance,  and  $N$ the size of the matrix $\mathcal{J}_{n}$.

Following the work in \cite{calarioster}, during the evaluation of the $\varphi-$ functions, the stopping criterion is
 $$10^{p} \cdot\Vert e_{m}\Vert< 1,$$ where $\,p$ being  the order of convergence of  the method ($p=2$ for the scheme EREM). In order to filter possible oscillations in the error estimate, 
the average on the last five  previous values  of the errors is used instead of  $ \Vert e_{m}\Vert$ in the  stopping condition. 
In the case of an unaccepted degree $m$, we increase the degree of the polynomial following relation  \eqref{newt}. When the degree $m$ for convergence is too large, the time step $\tau_{n}$ has to be split as described in \cite{calarioster}.
The algorithm with the function $\varphi_{1}$ is given in \cite{parallelleja}.


For the case where the spectral of $\mathcal{J}_{n}$ is more spread along  the imaginary axis, as for example  in  some hyperbolic problems, the method has been upgraded in \cite{ostfree}.

The attractive computational features of the method are clear  in  the sense that there is no Krylov subspace to store or linear systems to solve,
 but a drawback is that the method is based on interpolation, which is  generally ill-conditioned.  A major drawback is  that the required degree of the polynomial grows with the norm of the matrix $ \tau_n \mathcal{J}_{n}$.

\section{Numerical Simulations}
In the two examples, we deal with temperatures between 0 and $100 ^{o} \,\text{C} $. The water thermal expansivity $\alpha_{f}$
 and its compressibility $\beta_{f}$ used is from \cite{finem}. These two state functions depend on both pressure and the temperature.
 As we are dealing with low-enthalpy reservoirs ($T< 150 ^{o} \,\text{C} $), some water  
 properties can be well approximated as a function of temperature only. The water density in $kg.m^{-3}$ and the fluid viscosity  in $kg.m^{-1}.s^{-1}$ ( see \cite{graf})
 used are given respectively  by
$$ \rho_{f}(T)=1000 \left( 1- \left(\dfrac{\left(T-3.9863\right)^{2}}{508929.2}\times \dfrac{T+288.9414}{T+68.12963}\right)\right), $$ and 
\begin{eqnarray}
 \mu(T) = \left\lbrace \begin{array}{l}
 1.787 \times 10^{-3}\, e^{\left((-0.03288 + 1.962 \times 10^{-4} \times T)\times T\right)}\;\;\;\; \text{for}\;\;\, 0^{o} \,\text{C} \leq T\leq 40^{o} \,\text{C} \\
 \newline\\
 10 ^{-3} \times \left(1 + 0.015512 \times (T-20)\right)^{-1.572}\;\;\;\; \text{for}\;\;\,4 0^{o} \,\text{C}   < T\leq 100^{o} \,\text{C} \\
 \newline\\
0.2414 \times 10^{\left(\dfrac{247.8}{T + 133.15}\right)}\times 10^{-4}\;\;\;\; \text{for}\;\;\,100^{o} \,\text{C}  < T\leq 300^{o} \,\text{C} 
                                      \end{array}\right.
\end{eqnarray}
The water heat capacity in $J/kg.^{o} \text{C}$ used is also function of the temperature only in the interval $(0, 100)$ and  given by (see \cite[pp.73]{BundschuhArriaga})
\begin{eqnarray}
 c_{p}(T)=-1.3320081 \times 10^{-4}\,  T^{3}+0.0328405\, T^{2}-1.9254125\, T + 4206.3640128.
\end{eqnarray}
The water thermal conductivity used is $k_{f}=0.6 W/(m.K)$. We take the heat transfer coefficient $h_{e}$ sufficiently large to reach the local equilibrium.   
All our tests were performed on a workstation with a
3 GHz Intel processor and 8 GB RAM. Our code was
implemented in Matlab 7.11. We also used part of
the codes in  \cite{mrst}  for the spatial discretization. The absolute tolerance in the 
Krylov subspace  technique, Leja point technique, Newton iterations  and all linear systems is  $tol=10^{-6}$.  The dimension in the Krylov technique used is $m=10$.
The initial pressure used is the steady state pressure with
water properties at the initial temperature. 

In the legends of all of our graphs we use the following notation
\begin{itemize}
  \item ``Implicittheta=1'' denotes results from the  theta  Euler  with $\theta=1$ in\eqref{spaceheat2} and \eqref{spaceheat3}.
  \item  ``Implicittheta=0.5'' denotes results from the  theta  Euler  with $\theta=1/2$ in\eqref{spaceheat2} and \eqref{spaceheat3}.
  \item ``EREMKLeja'' denotes results from EREM  scheme with Krylov subspace for matrix exponential in the pressure system \eqref{spaceheat3} and real fast Leja points 
 for matrix exponential in the temperature system \eqref{spaceheat2}.
  \item ``EREMKrylov'' denotes results from EREM scheme  with Krylov subspace for matrix exponential in \eqref{spaceheat2} and \eqref{spaceheat3}.
  \item ``ROSM(1)'' denotes results from the  the scheme ROSM with $\gamma=1$ in\eqref{spaceheat2} and \eqref{spaceheat3}.
  \item ``ROSM(1/2)'' denotes results from the  the scheme ROSM with $\gamma=1/2$ in\eqref{spaceheat2} and \eqref{spaceheat3}.
 \item ``ROS2'' denotes results from the  scheme ROS2(1) in\eqref{spaceheat2} and \eqref{spaceheat3}.
\item ``ROS3p'' denotes results from the  scheme ROS3p in\eqref{spaceheat2} and \eqref{spaceheat3}.
\end{itemize}
 We use different constant time steps with the goal to study the convergence of the temperature and pressure equation at the final  time along with
the efficiency of numerical schemes.
The reference solutions used in the calculation of the errors are the numerical solutions with the
time step size equal to the half of the lower time step in the graphs.

\subsection{Heterogeneous  3D Geothermal Reservoir Simulation}
We consider a heterogeneous reservoir described by the domain
$\Omega=\left[0,1\right] \times \left[0,1\right] \times \left[0,0.1\right],$ where all
distances are in km. The half upper part of the
reservoir is less permeable than the lower part. The
injection point is located at the position $(1,1,z), \, z \in \left[0,1\right],$ injecting with the rate 
$q_{i}=1.04 m^{3}/s$, and the production is at $(0,0,z), \, z \in \left[0,1\right],$ with rate $q_{p}=-0.104 m^{3}/s$ with the lowest 
pressure at the point $(0,0,0)$.

Homogeneous Neumann boundary conditions are
applied for both the pressure equation, given by the
mass conservation law, and the energy equation. The water temperature at the injection well is $10 ^{o} \,\text{C}$. The
upper half of the reservoir has rock properties: permeability $K=10^{-2}$ Darcy, 
porosity $ \phi=20\%$, \,$\rho_{s}=2800 \,kg/m^{3}$, \,$c_{ps}=850 \,J/(kg.K)$,\,$k_{s}=2W/(m.K)$ while rock properties
of the lower half are: permeability  $K=10^{-1}$ Darcy, 
porosity $ \phi=40\%$, \,$\rho_{s}=3000 \,kg/m^{3}$, \,$c_{ps}=1000 \,J/(Kg.K)$,\,$k_{s}=3W/(m.K)$.
In  the two part the bulk vertical compressibility  $\alpha_{b}=10^{-7} \text{Pa}^{-1}$.  

 We  use a structured parallelepiped mesh. The size of the system is  $ 125000 \times 125000$ for the ODEs \eqref{spaceheat3} and  $ 250000\times 250000$ 
for the ODEs \eqref{spaceheat2}. The initial temperature at $z=0$ is   $60^{o} \,\text{C}$ and the temperature increases $ 3^{o} \,\text{C}$ at every 10 m.
 
\begin{figure}
 \subfigure[]{
   \label{FIG01a}
   \includegraphics[width=0.42\textwidth]{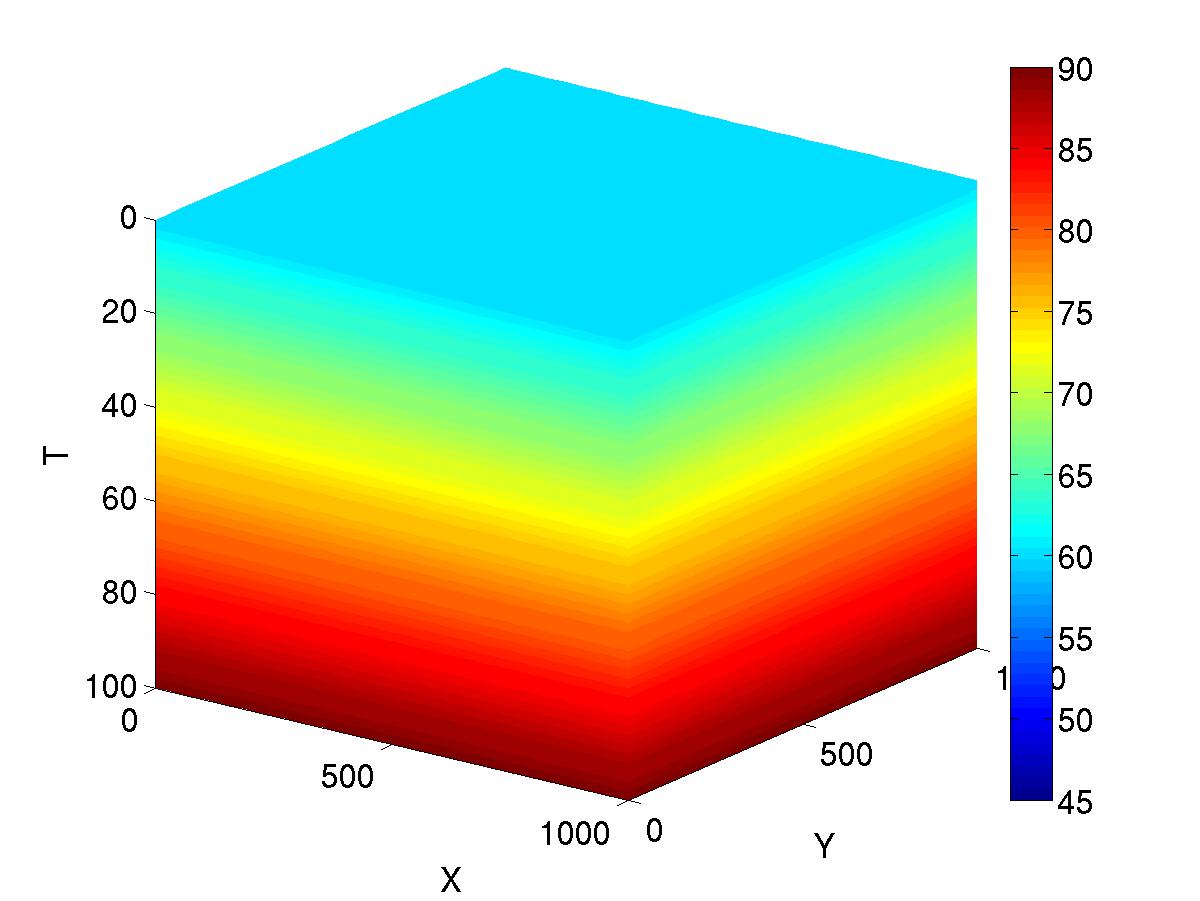}}
   \hskip 0.01\textwidth
   \subfigure[]{
   \label{FIG01b}
   \includegraphics[width=0.42\textwidth]{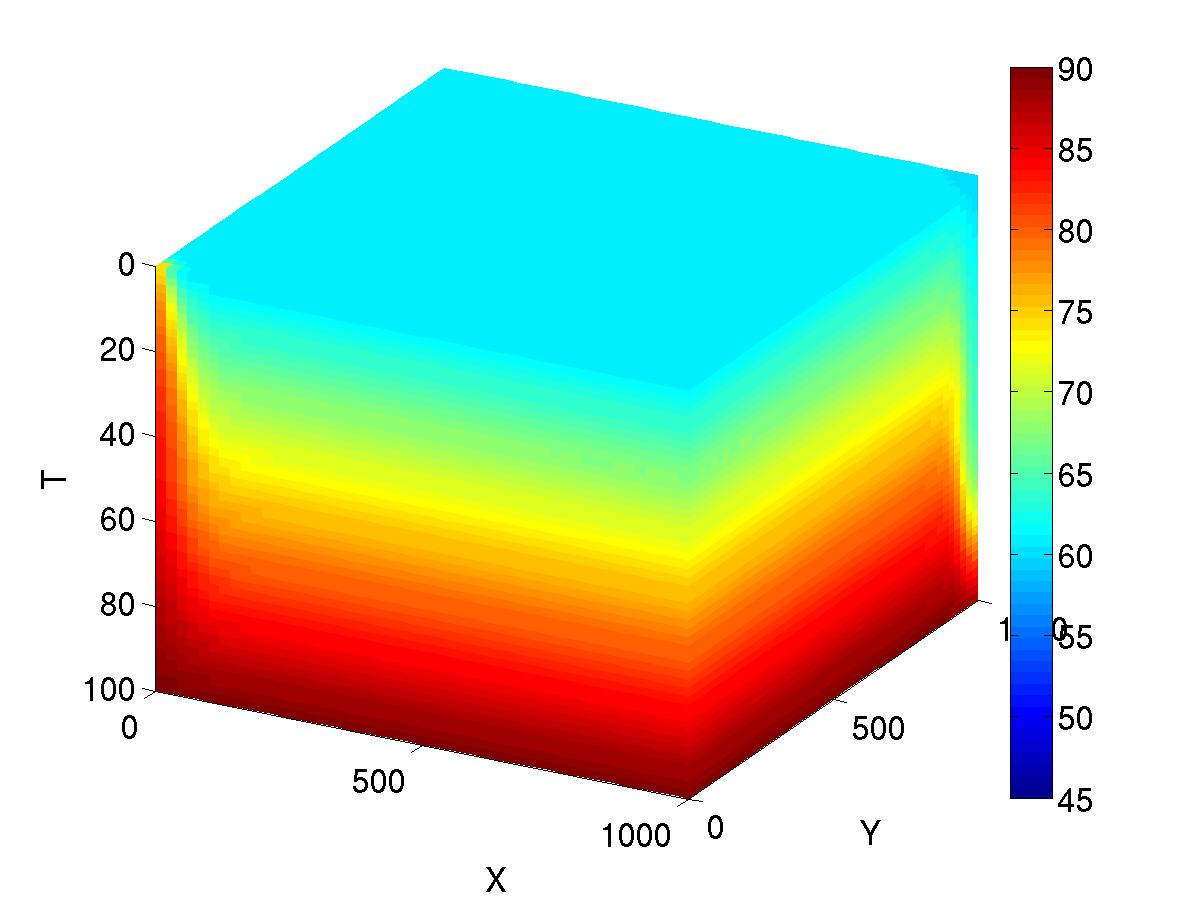}}
   \hskip 0.01\textwidth
   \subfigure[]{
   \label{FIG01c}
   \includegraphics[width=0.42\textwidth]{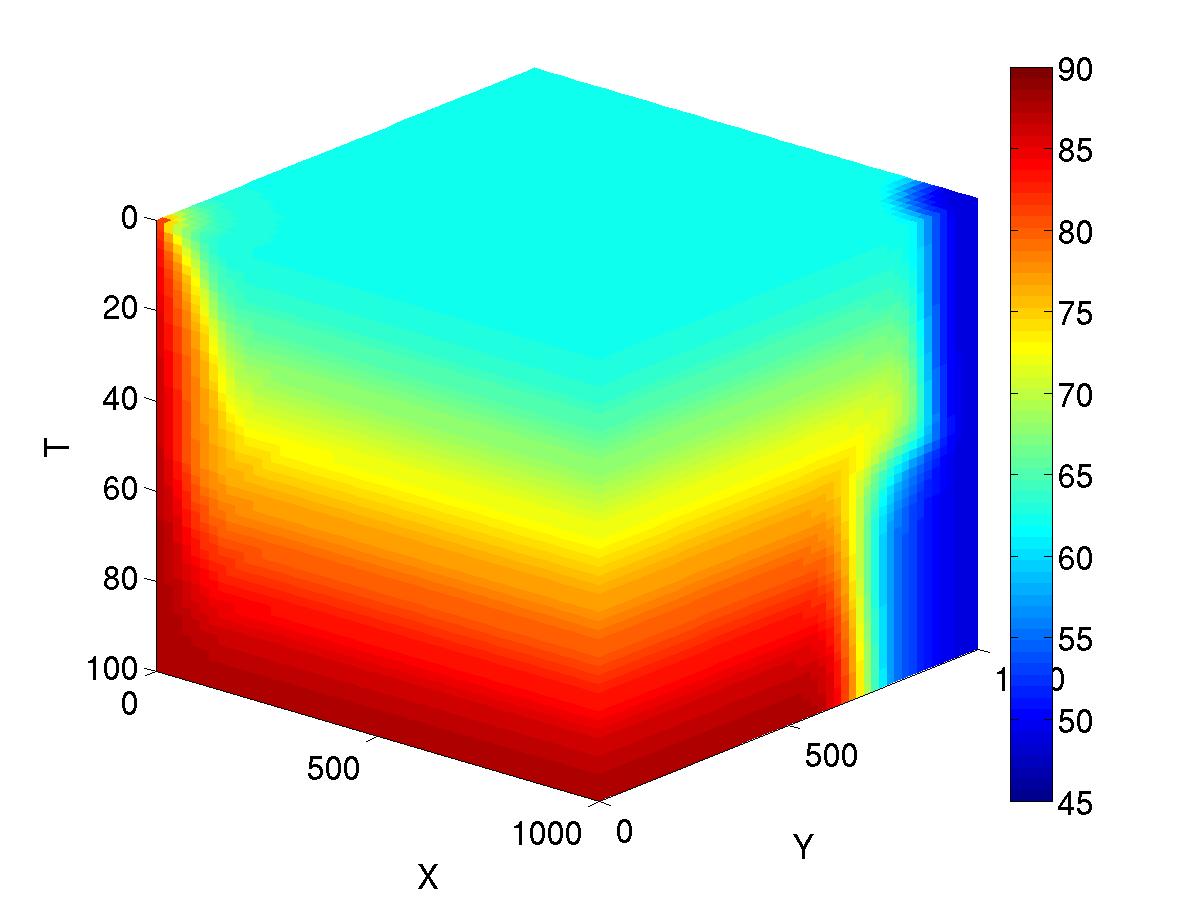}}
   \hskip 0.01\textwidth
   \subfigure[]{
   \label{FIG01d}
   \includegraphics[width=0.42\textwidth]{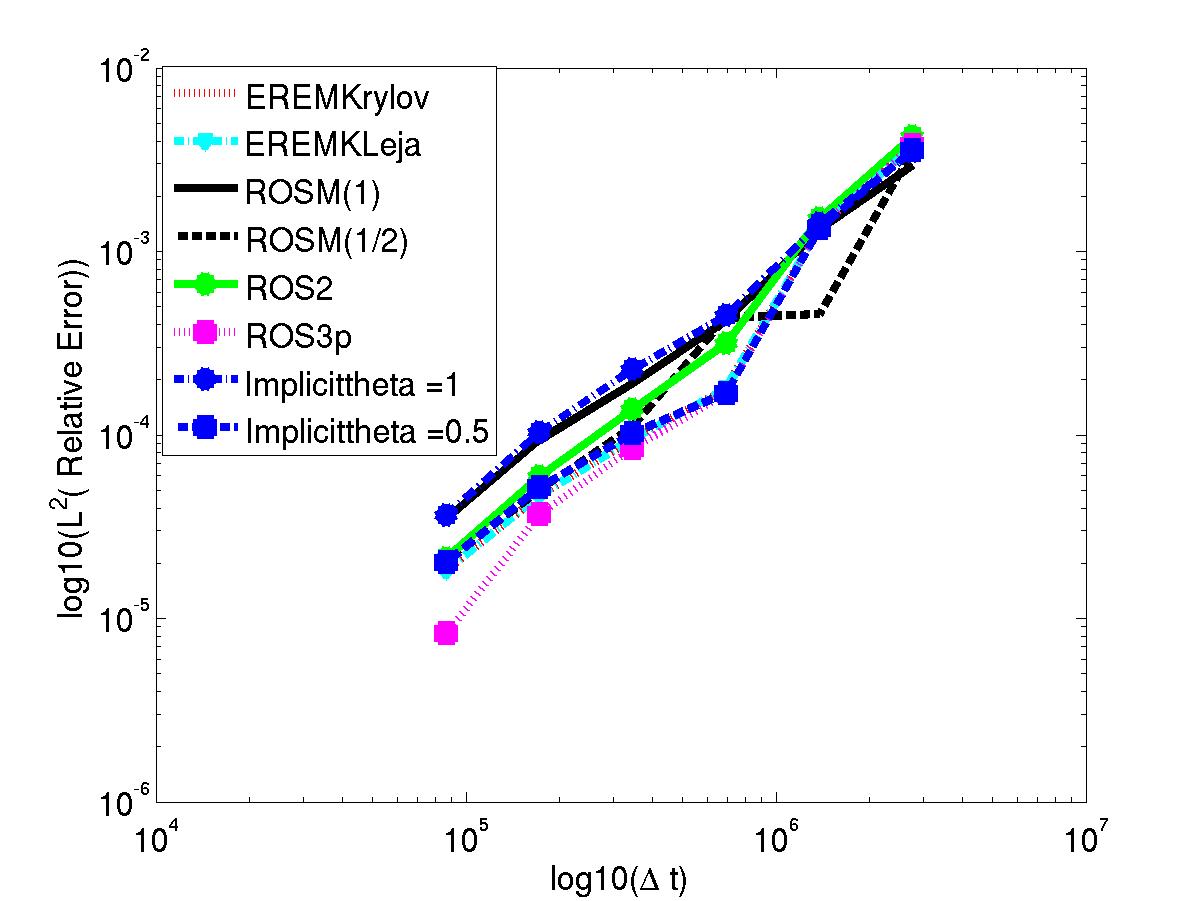}}
    \subfigure[]{
   \label{FIG01e}
   \includegraphics[width=0.43\textwidth]{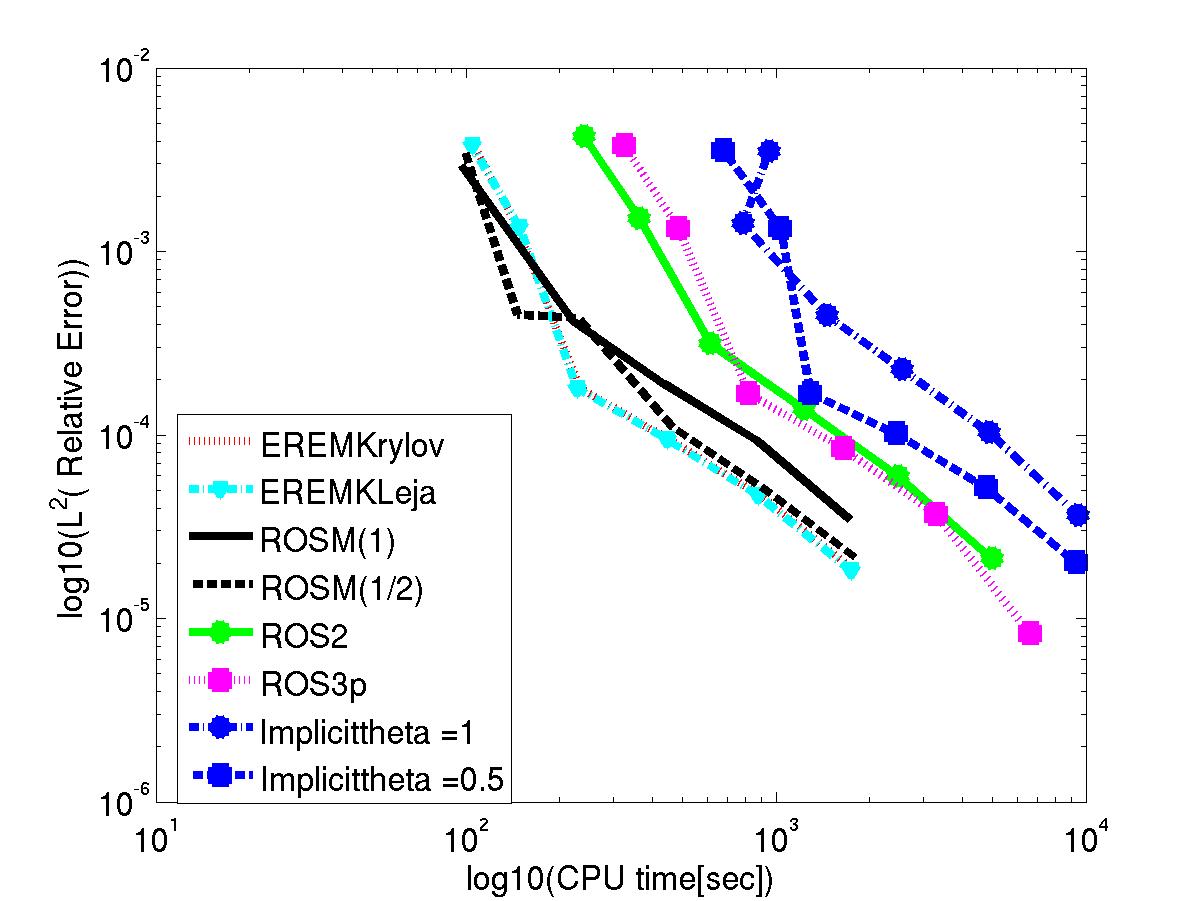}}
   \hskip 0.01\textwidth
   \subfigure[]{
   \label{FIG01f}
   \includegraphics[width=0.43\textwidth]{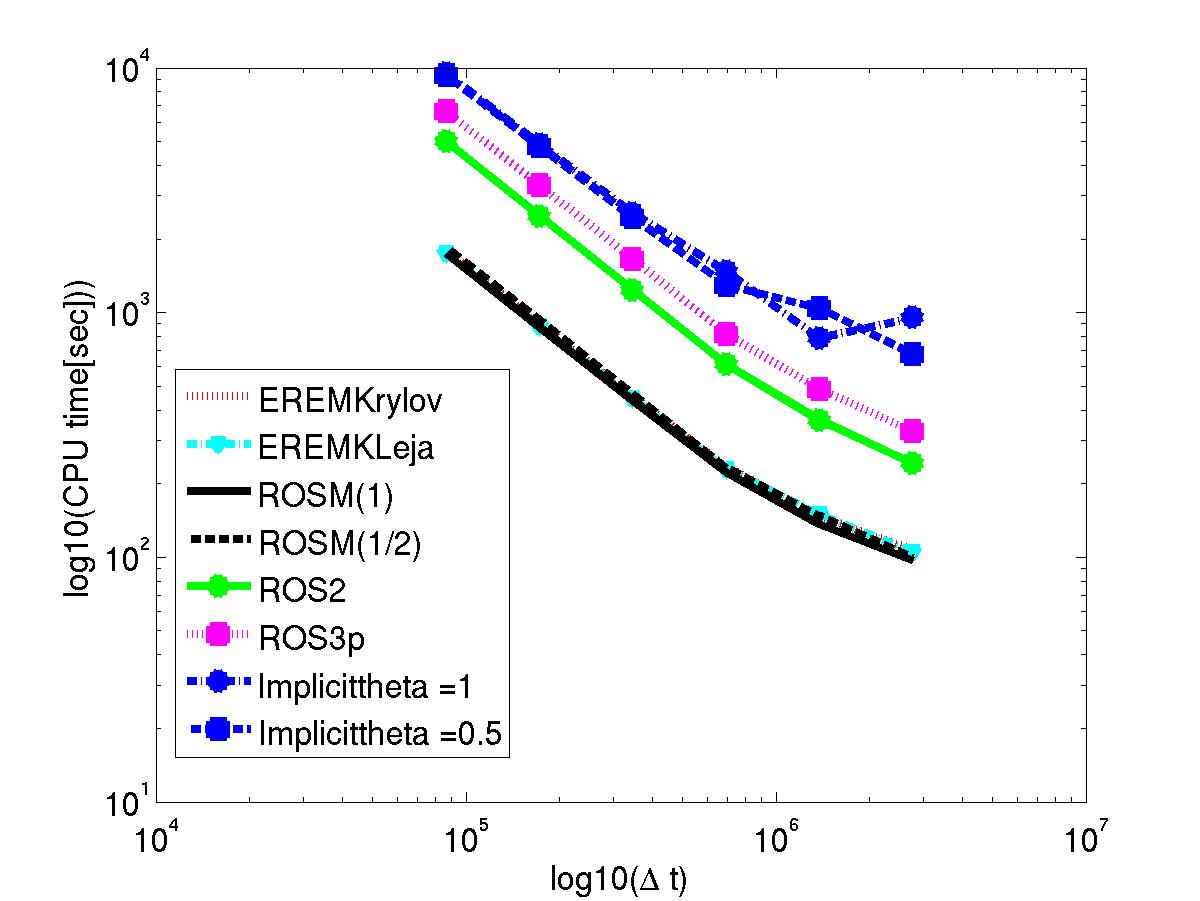}}
 \caption{(a) Initial temperature field, (b) temperature field at $\tau= 40$ days,  
(c) temperature field at  $\tau= 1000$ days, (d) the $L^{2}$ errors of the temperature as a function of time step size 
at the final time $\tau= 40$ days, (e) the corresponding $L^{2}$ errors as a function  of CPU time, and (f) CPU time as a function of time step size.
 }
 \label{FIG01}
\end{figure}

The initial temperature field is presented in \figref{FIG01a}, the temperature field at $\tau=40$ days is shown in \figref{FIG01b} 
while \figref{FIG01c}   shows the temperature at $\tau=1000$ days. We can  observe that the cold water decreases reservoir 
temperature at injection well and that the temperature at the production well increases.

 \figref{FIG01d} shows how the temperature errors at the final time $\tau=40$ days decrease with time step size.
 From this figure we can observe
 that the schemes with the same order of convergence in time have almost the same errors with our sequential approach. We can also observe that for large time steps, the 
errors are almost the same for all the schemes.
 The implicit $\theta$-Euler method with $\theta=1$ and the ROSM(1) 
are both of order $1.25$ in time. This order may decrease  up to 1 for  less smooth solutions \cite{lie}
 or relatively small time steps as for 
simple problems these schemes are order 1.
We can observe that the EREM scheme and the implicit $\theta$-Euler method with $\theta=1/2$ 
are slightly more accurate than the ROS2 scheme.
 EREM scheme and the implicit $\theta$-Euler method with $\theta=1/2$ are $1.55$ in time, 
 the ROS2(1) scheme is order $1.52$ and the ROS3p scheme is order $1.75$. This order may decrease for 
less smooth solutions \cite{lie}. We can however observe that schemes with high orders in time 
for simple problems are affected  by order reduction in the sequential approach.

\figref{FIG01e} shows the relative $L^{2}$  temperature errors  as the function of the CPU time corresponding to \figref{FIG01d}.  We can observe the efficiency of the EREM scheme comparing the others schemes.
This figure also shows that the Exponential Rosenbrock-Euler Method and Rosenbrock- Type methods are very efficient as compared to the standard implicit $\theta$-Euler
methods.

\figref{FIG01f} shows the CPU time as a function of time step size corresponding to \figref{FIG01d} and \figref{FIG01e}. We can 
clearly observe that the Exponential Rosenbrock-Euler Method and Rosenbrock- Type methods are again very efficiency compared to standard implicit $\theta$-Euler methods.
From this figure we can observe that EREM scheme, ROSM(1) and ROSM(1/2) schemes are at least five time as efficient, the ROS2 is at least  
twice as efficient and the ROS3p at least one and half time as efficient as the standard methods. While these factors may be dependent on the particular implementation, and the availability of good nonlinear solvers for the nonlinear solves in the standard methods, we believe them to be representative.  

\subsection{2D Fractured Geothermal Reservoir Simulation}
We consider here a 2 D fractured reservoir \cite{torh,TPFAdfm}, with a quasi-structured triangular mesh in the domain
$\Omega=\left[0,100\right] \times \left[0,100\right] $, where all
distances are in m. The  matrix  properties are:  permeability  $10^{-2}$ Darcy, porosity  $\phi=20\%$,\, $c_{ps}=1000 \,J/(Kg.K)$,\,$\rho_{s}=2800 \,kg/m^{3}$, $k_{s}=2W/(m.K)$ and $\alpha_{b}=0$. 
The fractures have an aperture of $10^{-3}m$ and a permeability of $100 $ Darcy.  The
injection point is located at the position $(0,40)$,  with constant pressure $10$MPa, while the production point is located at the point  $(100,100)$, with constant pressure $10^5$Pa.

Homogeneous Neumann boundary conditions are
applied for both the pressure equation given by the
mass conservation law and the energy equation.
The size of the system is  $ 11460 \times 11460$ for the ODEs \eqref{spaceheat3} and  $ 22920\times 22920$ 
for the ODEs \eqref{spaceheat2}. The initial temperature is $80 ^{o} \,\text{C}$ while
the water temperature at the injection well is $15 ^{o} \,\text{C}$.

\begin{figure}
  \subfigure[]{
\label{FIG02a}
   \includegraphics[width=0.5\textwidth]{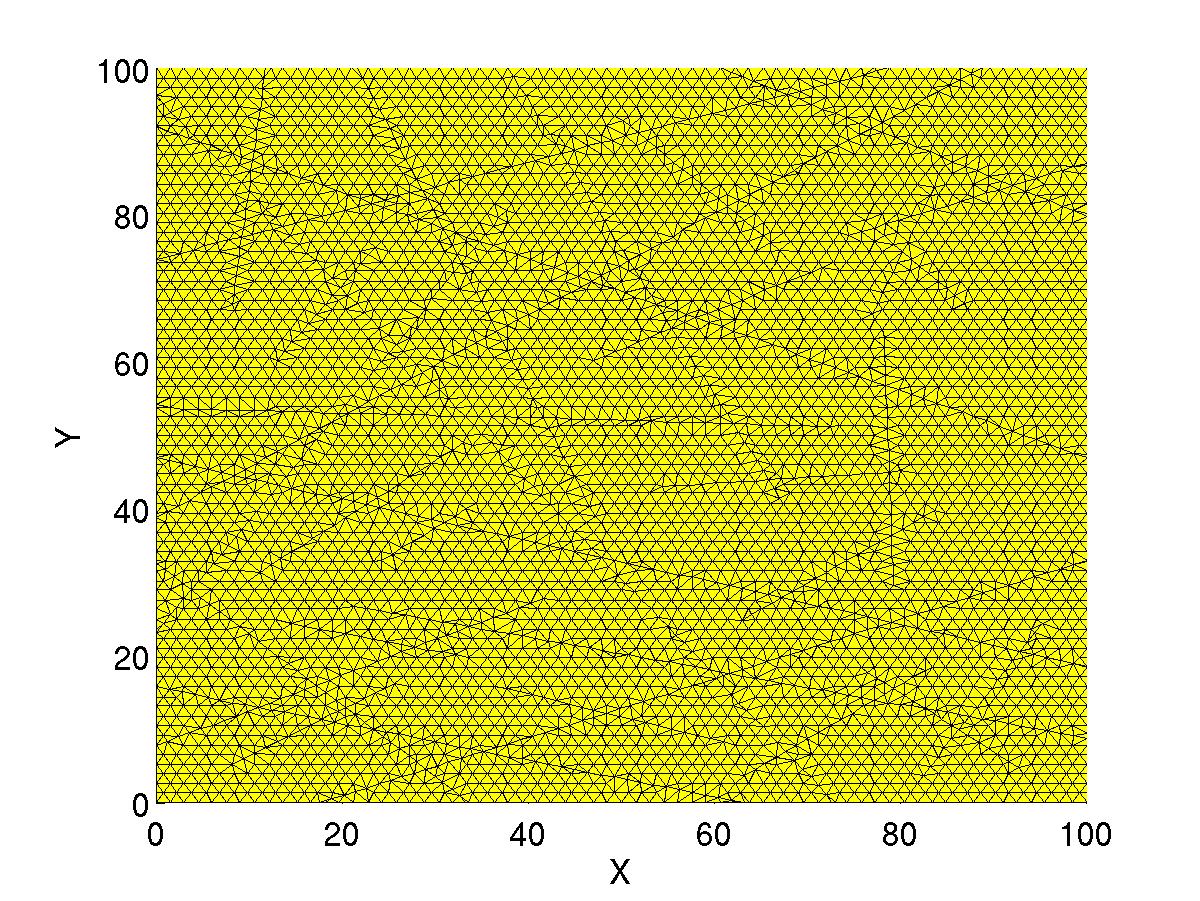}}
   \hskip 0.01\textwidth
   \subfigure[]{
   \label{FIG02b}
   \includegraphics[width=0.5\textwidth]{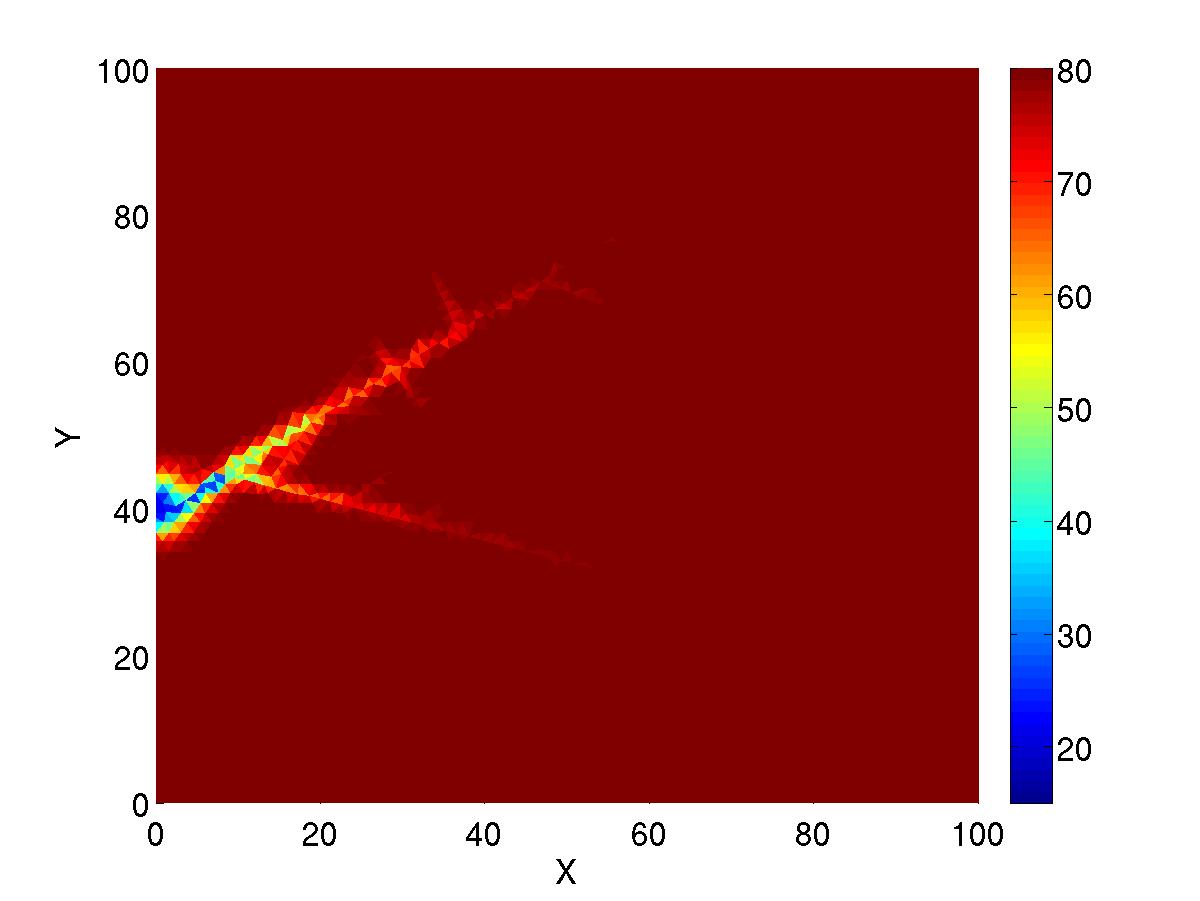}}
   \hskip 0.01\textwidth
\subfigure[]{
   \label{FIG02c}
   \includegraphics[width=0.5\textwidth]{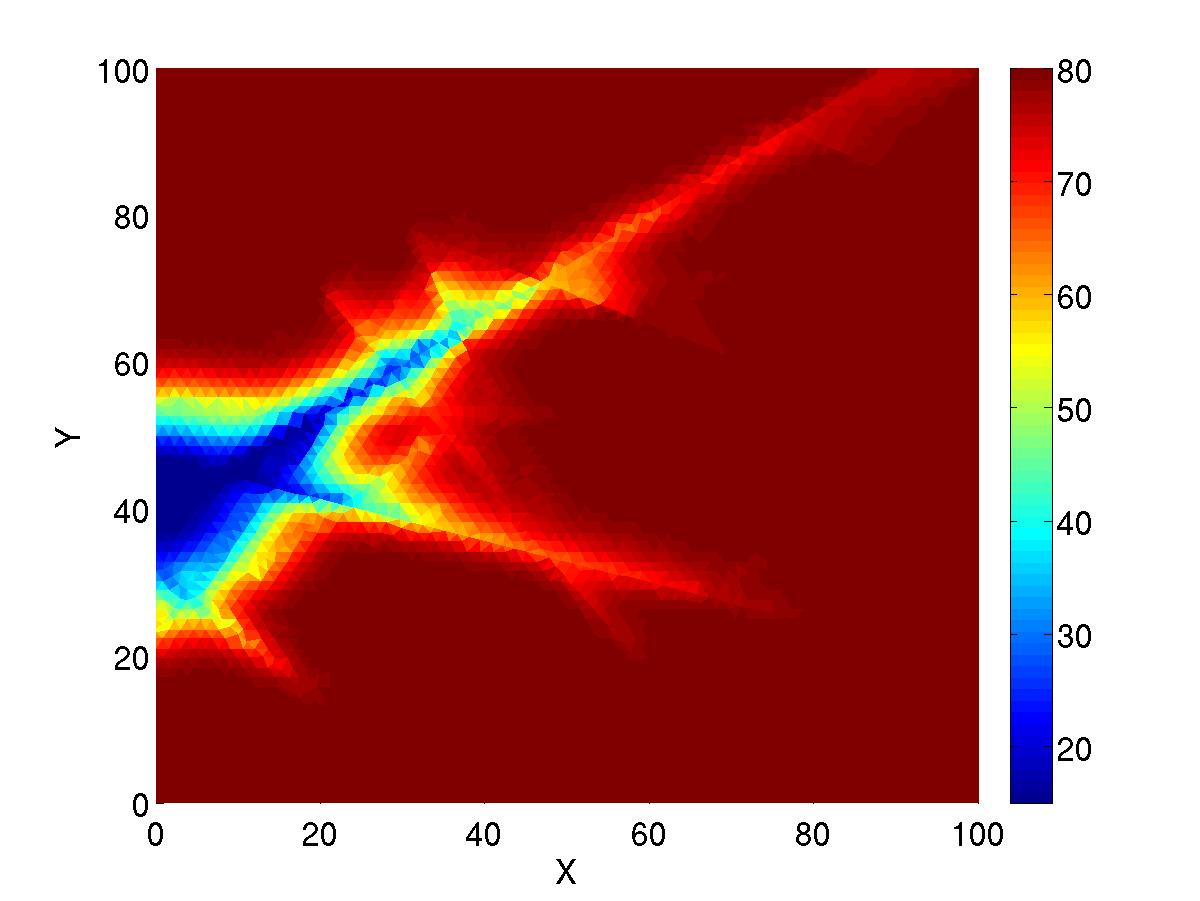}}
   \hskip 0.01\textwidth
\subfigure[]{
   \label{FIG02d}
   \includegraphics[width=0.5\textwidth]{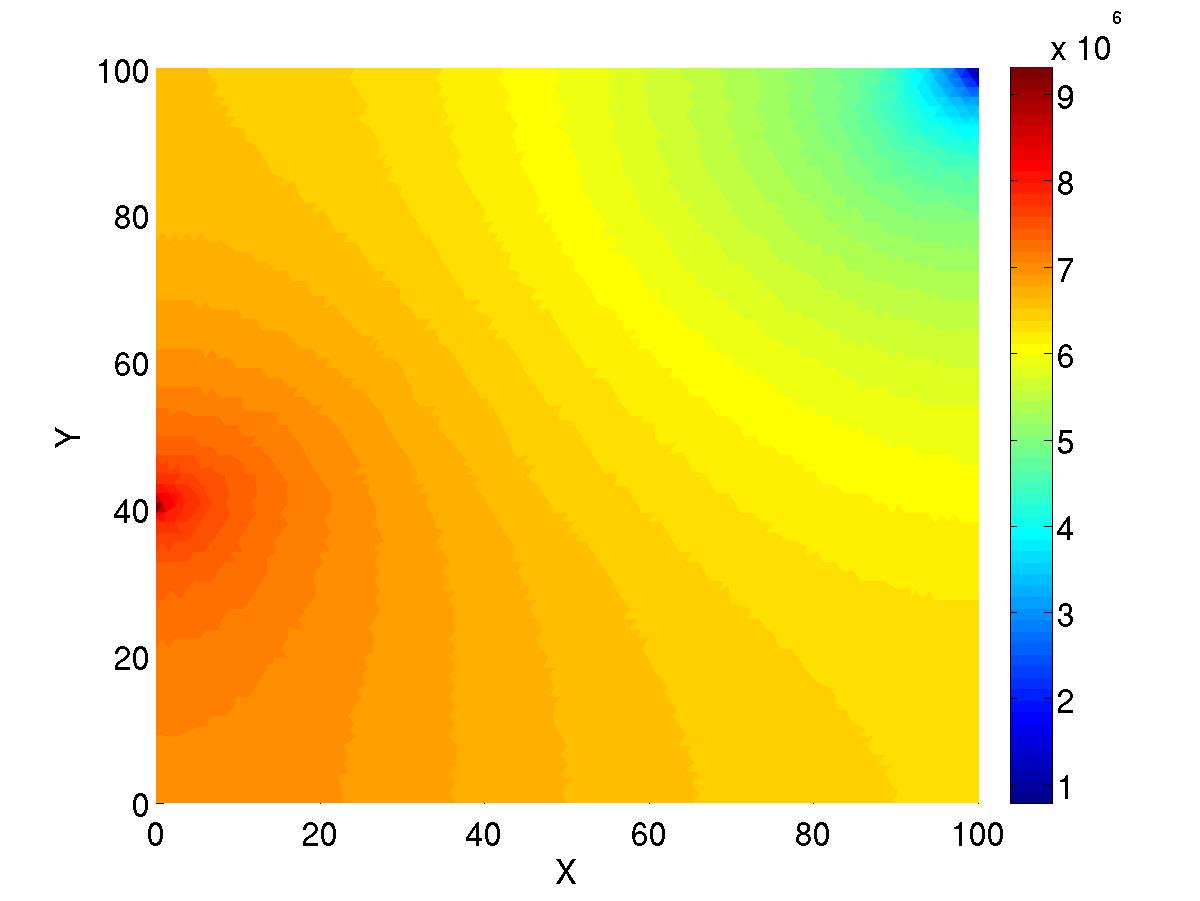}}
\caption{ (a) Fractured 2D grid, (b) temperature field at time $\tau=10$ days, (c) temperature field at 
the time $\tau=100$ days, and (d) the pressure field at time $\tau=10$ days. The injection well is at the point $(0,40)$ 
and the production well at  the point $(100,100)$.}
 \label{FIG02}
\end{figure}

The 2D grid with fractures is shown in \figref{FIG02a}, the temperature field at time $\tau=10$ 
 days  in \figref{FIG02b} and the temperature field at the time $\tau=100$ days in \figref{FIG02c}.  \figref{FIG02d} 
shows the pressure field  at time $\tau=10$ days.

\begin{figure}
 \subfigure[]{
   \label{FIG03a}
   \includegraphics[width=0.5\textwidth]{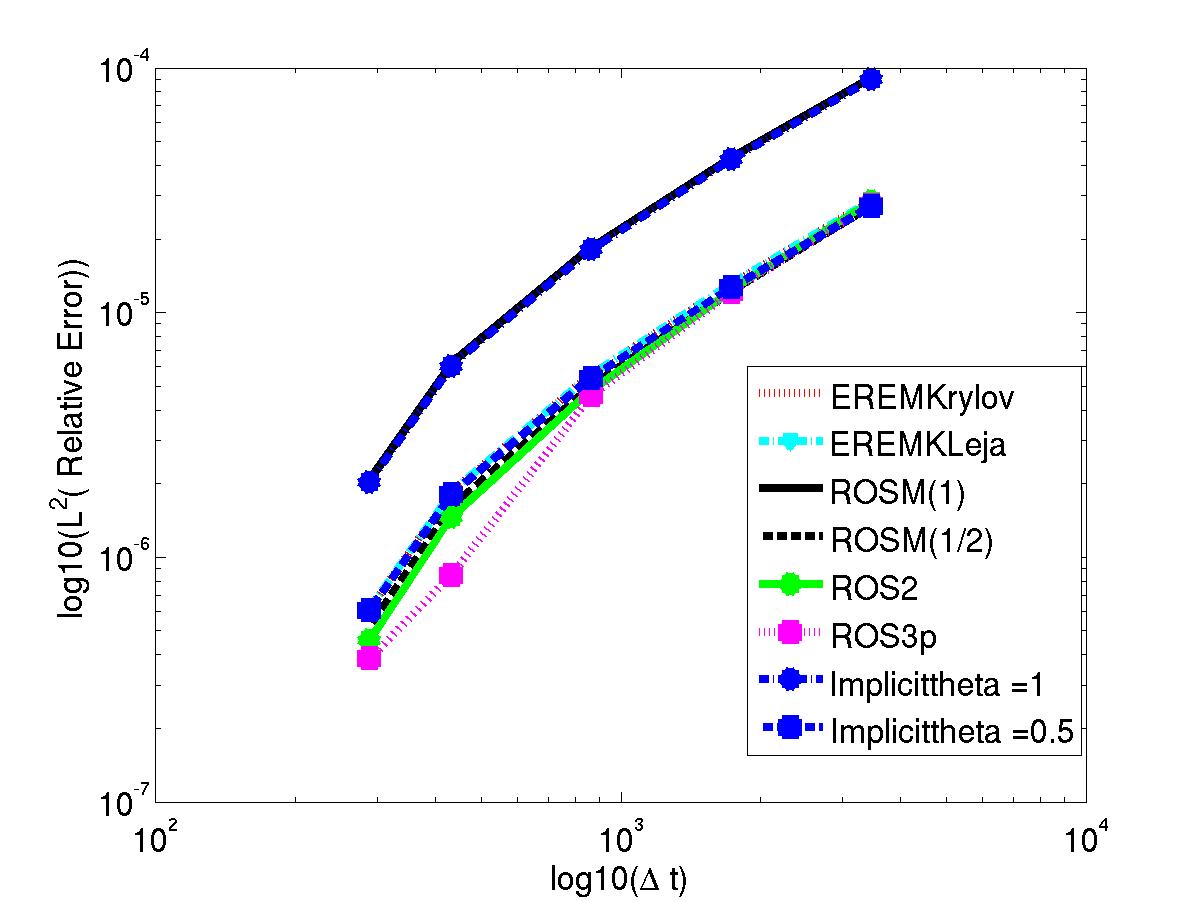}}
   \hskip 0.01\textwidth
   \subfigure[]{
   \label{FIG03b}
   \includegraphics[width=0.5\textwidth]{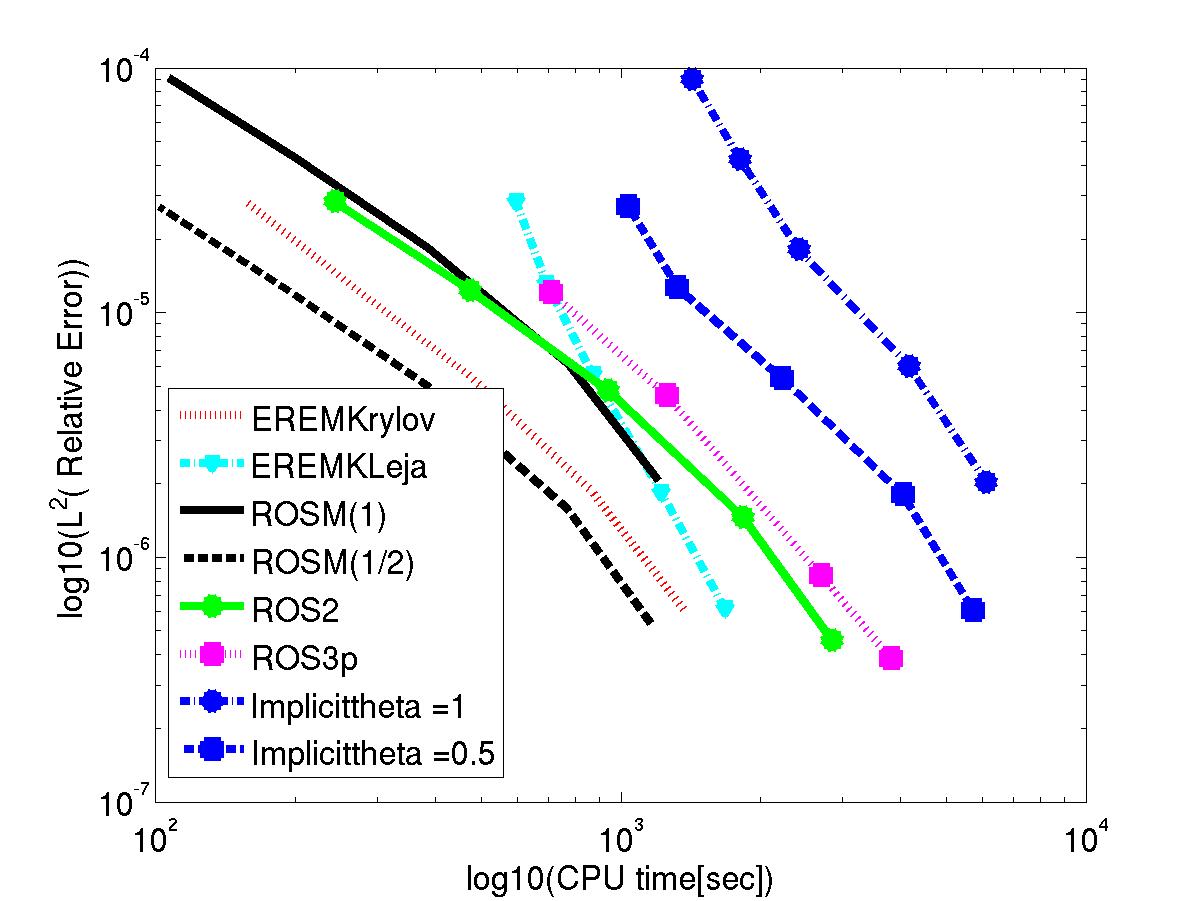}}
    \subfigure[]{
   \label{FIG03c}
   \includegraphics[width=0.5\textwidth]{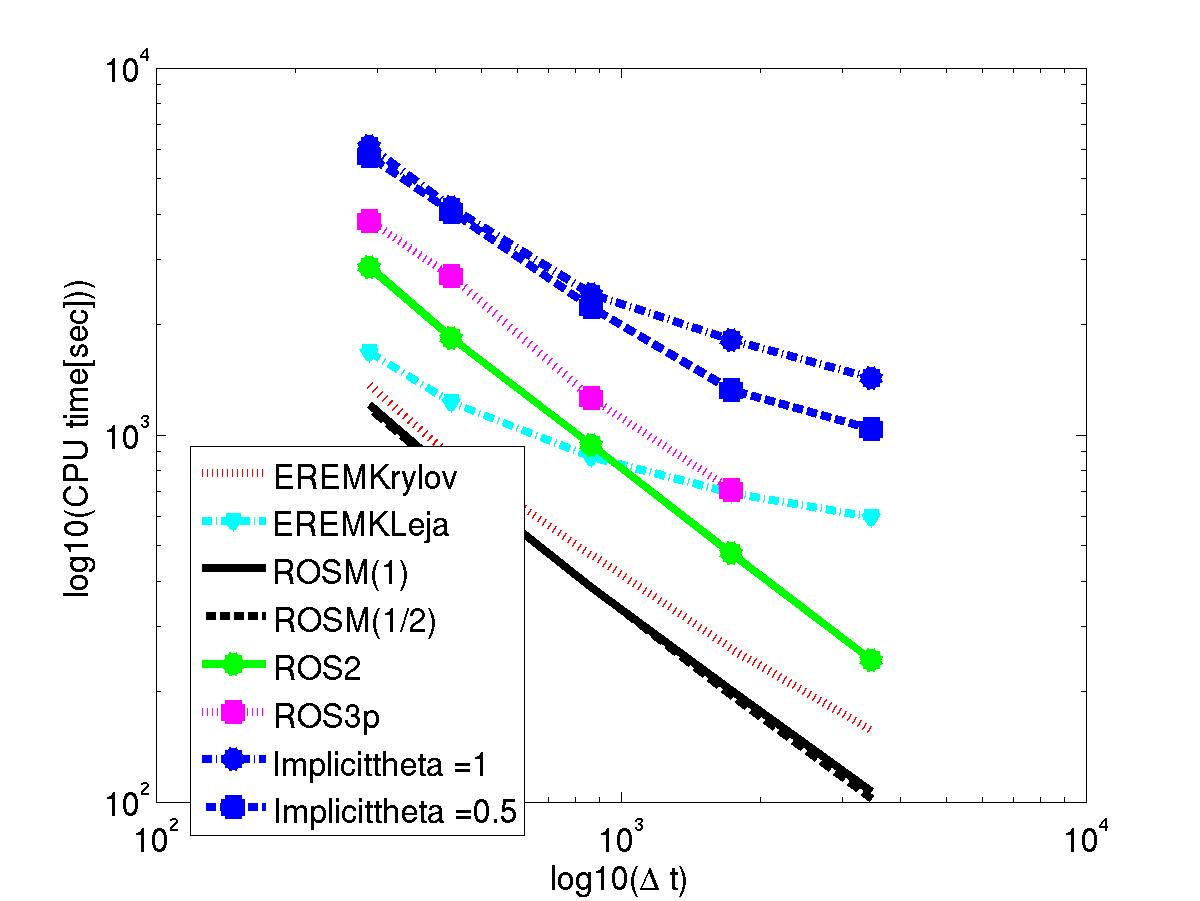}}
   \hskip 0.01\textwidth
   \subfigure[]{
   \label{FIG03d}
   \includegraphics[width=0.5\textwidth]{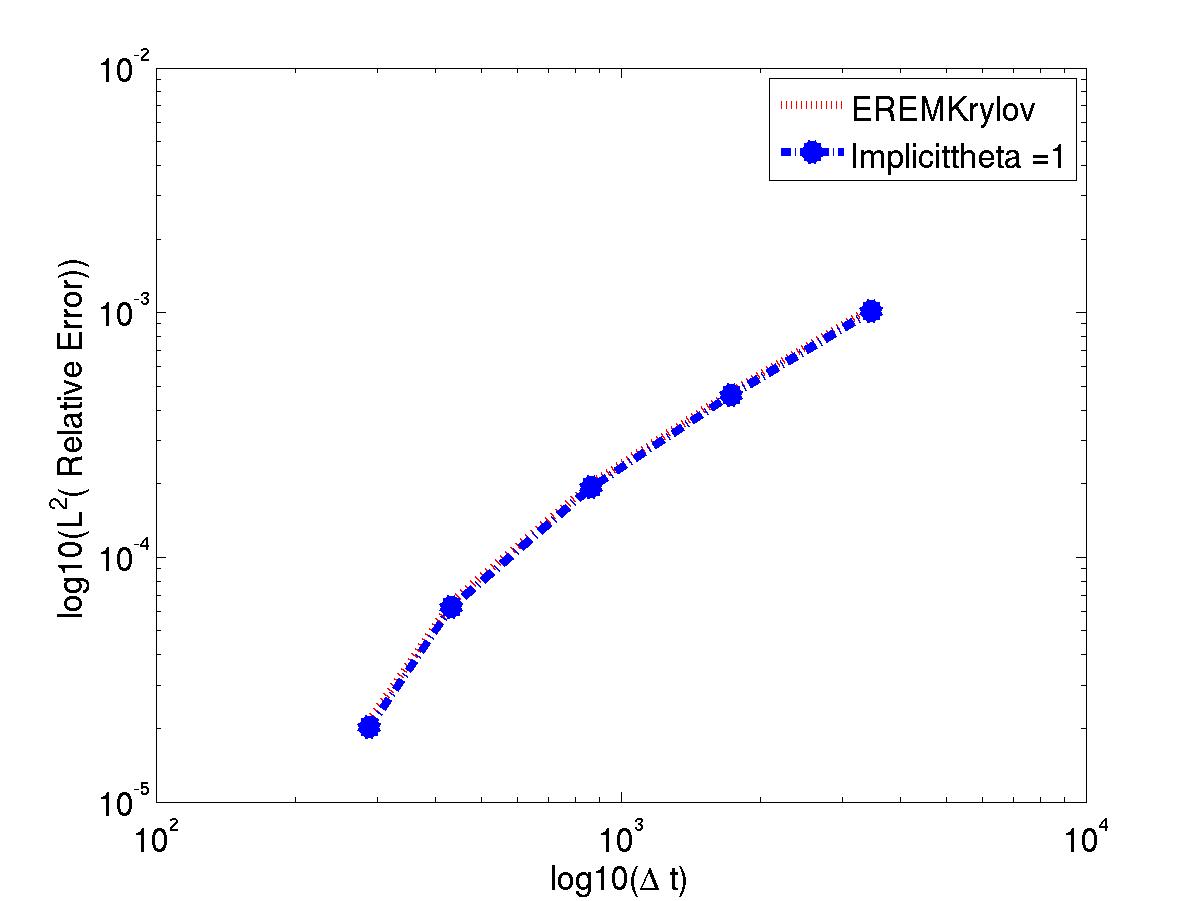}}
 \caption{
(a) The relative $L^{2}$ errors of the temperature as a function of time step size 
at the final time $\tau= 10$ days, (b) the corresponding $L^{2}$ errors as a function  
of CPU time, (c) CPU time as a function of time step size, and
(d) the relative $L^{2}$ errors of the pressure as a function of time step size 
at the final time $\tau= 10$ days, we plot only for two schemes as all the methods give the same errors.
}
 \label{FIG03}
\end{figure}

\figref{FIG03a} shows the time convergence of the all 
the schemes as the  temperature errors decrease with time step size at the final time $\tau=10$ days.
 From this figure we can observe again
 that the schemes with the same order of convergence in time have almost  the same errors. 
 The implicit $\theta$-Euler method with $\theta=1$ and the ROSM(1) 
are order $1.4$ in time.
 The EREM scheme, the $\theta$-Euler method with $\theta=1/2$ and the ROS2(1) scheme are  order $1.50$ in time while the ROS3p scheme is order $1.75$.
 These orders may increase for relatively small time steps. Again we can observe that schemes with high orders in time 
for simple problems are affected  by order reduction in the sequential approach.

\figref{FIG03b} shows the relative $L^{2}$  temperature errors  as function of the CPU time corresponding to \figref{FIG03a}. 
As in the first example, we can observe the efficiency of the schemes ROSM(1/2) and EREM  with the Krylov technique compared to the others schemes.

\figref{FIG03c} shows the CPU time as a function of time step size, corresponding to figure \figref{FIG03a} and \figref{FIG03b}. Again we observe that EREM scheme, ROSM(1) and ROSM(1/2) schemes are  almost four times as efficient as the standard implicit methods, while the ROS2 and the ROS3p are almost twice as efficient as the standard implicit methods. From these examples, we expect the efficiency gain to increase with the size of the problem.

The relative $L^{2}$ pressure errors are almost the same for all the numerical schemes. 
We therefore plot only the errors for two numerical schemes with order 1 and 2  in time respectively.
 \figref{FIG03a} shows the pressure errors at time $\tau=10$ days  as a function of time step size for the EREM  scheme and $\theta$-Euler method with $\theta=1$. 
We can observe the convergence of those schemes while solving the pressure equation.
The order of convergence in time is $1.5$ and may decrease according to
\cite{lie} for rough solutions (less smooth solutions).

\section{Conclusion}
We have proposed a novel approach for simulation of
geothermal processes in heterogeneous porous media.
This approach decouples the mass conservation equation from
the energy equation and solves each stiff ODEs from
space discretization sequentially using the Exponential Rosenbrock-Euler method and Rosenbrock-type methods for the time integration.

Numerical simulations in 2D and 3D show that using the Krylov
subspace technique and Real Leja points technique in the computation of the
exponential functions $\varphi_{i}$ in the Exponential Rosenbrock-Euler method, and  the 
Matlab function bicgstab with ILU(0) preconditioners with no fill-in for solving all linear systems appearing in the Rosenbrock-type methods and the implicit Euler theta methods, makes our
approach more efficient compared to the sequential standard implicit Euler methods.

\section*{ACKNOWLEDGEMENTS}
We thank Tor Harald Sandve for sharing with us the 2D fractured grid and its discretization.
 This work was funded by the Research Council of
Norway (grant number 190761/S60).

 
\end{document}